\documentclass[12pt]{article}

\usepackage{times}

\usepackage{enumitem}
\usepackage[utf8]{inputenc}
\usepackage{commath}
\usepackage{amsthm, amscd}
\usepackage{amsmath}
\usepackage{amssymb}
\usepackage{cite}
\usepackage{setspace}
\usepackage{ stmaryrd }
\usepackage{tikz-cd}

\usepackage{tocloft}
\usepackage{sectsty}
\usepackage{xparse}

\newcommand*\tasklabelformat[1]{#1)}

\usepackage{titlesec}

\usepackage{tasks}[newest]
\settasks{
  label = \alph* ,
  label-format = \tasklabelformat ,
  label-width  = 12pt
}

\usepackage{bigints}
\usepackage{comment}
\usepackage{mathtools}
\usepackage{mathrsfs}
\usepackage{fancyhdr}

\allowdisplaybreaks[4]

\numberwithin{equation}{section}
\usepackage{etoolbox}
\patchcmd{\thebibliography}
  {\settowidth}
  {\setlength{\itemsep}{0pt plus -10pt}\settowidth}
  {}{}
\apptocmd{\thebibliography}
  {
  }
  {}{}
\makeatletter
\newtheorem*{rep@theorem}{\rep@title}
\newcommand{\newreptheorem}[2]{%
\newenvironment{rep#1}[1]{%
 \def\rep@title{#2 \ref{##1}}%
 \begin{rep@theorem}}%
 {\end{rep@theorem}}}
\makeatother

\theoremstyle{theorem}

\newreptheorem{theorem}{Theorem}
\newtheorem{thm}{Theorem}[section]
\newtheorem*{thm*}{Theorem}
\theoremstyle{definition}
\newtheorem{prop}[thm]{Proposition}
\newtheorem*{prop*}{Proposition}

\newtheorem{lem}[thm]{Lemma}
\newtheorem{cor}[thm]{Corollary}
\newtheorem*{cor*}{Corollary}
\theoremstyle{remark}

\newtheorem{rem}[thm]{Remark}

\title{\vspace*{-1.5cm}
On the Harder-Narasimhan slopes of direct images
}

\author
{Siarhei Finski
}

\date{}

\usepackage[%
    left=1in,%
    right=1in,%
    top=1.1in,%
    bottom=0.8in,%
    paperheight=11in,%
    paperwidth=8.5in%
]{geometry}

\newcommand{\sym}{{\rm{Sym}}}

\newcommand{\comp}{\mathbb{C}}
\newcommand{\real}{\mathbb{R}}

\newcommand{\nat}{\mathbb{N}}
\newcommand{\integ}{\mathbb{Z}}

\newcommand{\rk}[1]{{\rm{rk}} ( #1 )}

\makeatletter
\DeclareFontFamily{OMX}{MnSymbolE}{}
\DeclareSymbolFont{MnLargeSymbols}{OMX}{MnSymbolE}{m}{n}
\SetSymbolFont{MnLargeSymbols}{bold}{OMX}{MnSymbolE}{b}{n}
\DeclareFontShape{OMX}{MnSymbolE}{m}{n}{
    <-6>  MnSymbolE5
   <6-7>  MnSymbolE6
   <7-8>  MnSymbolE7
   <8-9>  MnSymbolE8
   <9-10> MnSymbolE9
  <10-12> MnSymbolE10
  <12->   MnSymbolE12
}{}
\DeclareFontShape{OMX}{MnSymbolE}{b}{n}{
    <-6>  MnSymbolE-Bold5
   <6-7>  MnSymbolE-Bold6
   <7-8>  MnSymbolE-Bold7
   <8-9>  MnSymbolE-Bold8
   <9-10> MnSymbolE-Bold9
  <10-12> MnSymbolE-Bold10
  <12->   MnSymbolE-Bold12
}{}

\let\llangle\@undefined
\let\rrangle\@undefined
\DeclareMathDelimiter{\llangle}{\mathopen}%
                     {MnLargeSymbols}{'164}{MnLargeSymbols}{'164}
\DeclareMathDelimiter{\rrangle}{\mathclose}%
                     {MnLargeSymbols}{'171}{MnLargeSymbols}{'171}
                     
\DeclareMathOperator*{\esssup}{ess\,sup}
\DeclareMathOperator*{\essinf}{ess\,inf}
\makeatother

\newenvironment{sciabstract}{}

\cftsetindents{section}{0em}{1.7em}

\sectionfont{\large}

\titlespacing*{\section}
{0pt}{5pt}{5pt}

\begin{document}

\maketitle

\begin{sciabstract}
  \textbf{Abstract.}
  For a flat family of complex projective varieties, we study the asymptotic distribution of the Harder-Narasimhan slopes of the direct images on the base of high tensor powers of a line bundle on the total space. 
  We establish a theorem of Mehta-Ramanathan type, showing that this asymptotic distribution can be recovered from its analogues associated with the base changes of the family over general curves.
  \par 
	As an application, we show that the multiplicity of the line bundle on the family can be estimated in terms of the Harder-Narasimhan slopes of the associated direct images and the multiplicity of the general fiber.
\end{sciabstract}

\pagestyle{fancy}
\lhead{}
\chead{On the Harder-Narasimhan slopes of direct images}
\rhead{\thepage}
\cfoot{}


\newcommand{\Addresses}{{
  \bigskip
  \footnotesize
  \noindent \textsc{Siarhei Finski, Centre de Mathématiques Laurent Schwartz (CMLS), CNRS, École polytechnique, Institut Polytechnique de Paris, Palaiseau, France}\par\nopagebreak
  \noindent  \textit{E-mail }: \texttt{siarhei.finski@polytechnique.edu} / \texttt{finski.siarhei@gmail.com}.
}} 

\vspace*{0.25cm}

\par\noindent\rule{1.25em}{0.4pt} \textbf{Table of contents} \hrulefill

\vspace*{-1.5cm}

\tableofcontents

\vspace*{-0.2cm}

\noindent \hrulefill


\section{Introduction}\label{sect_intro}
	Consider a surjective flat holomorphic map $\pi : X \to B$ with irreducible general fibers between a complex projective irreducible variety $X$ and a projective manifold $B$ of dimensions $n + m$ and $m$ respectively, $n, m \in \nat^*$.
	Let $L$ be a holomorphic line bundle over $X$.
	We fix a \textit{multipolarization} $[\omega_B]$, which is a collection of Kähler classes $[\omega_{B, i}] \in H^{1, 1}(B)$, $i = 1, \ldots, m - 1$.
	In this paper, we study the Harder-Narasimhan $[\omega_B]$-slopes of direct image sheaves $\mathscr{E}_k := R^0 \pi_* L^k$ for big powers $k \in \nat$.
	\par 
	The study of Harder-Narasimhan slopes on direct images has recently received a lot of attention, especially when the base of the family is a curve and the line bundle $L$ is relatively ample.
	\par 
	Chen studied its applications in Arakelov geometry in \cite{ChenHNolyg} and algebraic geometry in \cite{ChenVolume}.
	It was applied in the study of the moduli space of K-stable Fano varieties by Codogni-Patakfalvi \cite{CodPatak}, Xu-Zhuang \cite{XuZhuPosCM} and more recently by Hattori \cite{HattoriCM}.
	Related ideas have been used in diophantine approximations, see Faltings-Wüstholz \cite{FaltingsWusth}, cf. Grieve \cite{GrieveDioph}.
	\par 
	The goal of the present paper is twofold.
	First, we show that several results previously known when the line bundle is relatively ample and the base of the family is a curve continue to hold in general.
	Second, we prove that the Harder-Narasimhan slopes exhibit a good behavior under pull-backs to general curves, giving an asymptotic analogue of the theorem of Mehta-Ramanathan.
	\par 
	The extension of the theory for more general bases is very natural, especially when considered in parallel with the recent developments concerning the Wess-Zumino-Witten equation from \cite{FinHNII} and \cite{FinHYM}.
	Indeed, the latter equation can be seen as a generalization of the Hermite-Einstein equation, cf. \cite[\S 1]{FinHYM}, and one of the main results of \cite{FinHYM} shows that the variability of the Harder-Narasimhan slopes on direct images provides the full set of obstructions to the existence of approximate solutions to the Wess-Zumino-Witten equation, in exactly the same way as the variability of the Harder-Narasimhan slopes of vector bundles provides the full set of obstructions to the existence of approximate solutions to the Hermite-Einstein equation, cf. \cite[\S 6]{KobaVB}.
	It is hence as natural to consider the measure associated with the Harder-Narasimhan slopes over a higher-dimensional base as it is to consider the Hermite-Einstein equation on a higher-dimensional manifold rather than a curve.
	\par 
	Moreover, in \cite{FinHNII} and \cite[\S 1]{FinHYM}, the author observed that the Wess-Zumino-Witten equation can alternatively be seen as the asymptotic version of the Hermite-Einstein equation, which allowed to develop a fibered analogue of Atiyah-Bott theory \cite{AtiyahBott}, identifying the infima of the associated Yang-Mills functionals with the quantities expressible in terms of the limiting Harder-Narasimhan measures, see \cite[Theorem 1.1]{FinHYM} for a precise statement.
	Both articles, \cite{FinHNII} and \cite{FinHYM}, appeared after the current one, and the results of the current article were used there.
	\par 
	Our motivation for developing the corresponding asymptotic Mehta-Ramanathan theorem also arose from its connection with the Wess-Zumino-Witten equation.
	Indeed, as the Wess-Zumino-Witten equation can be viewed as an asymptotic version of the Hermite-Einstein equation, it is natural to expect that the asymptotic Mehta-Ramanathan theorem could play a role analogous to that played by the classical Mehta-Ramanathan theorem in the Hermite-Einstein setting, used in the proof of Donaldson \cite{DonaldASD} of the Kobayashi-Hitchin correspondence.
	Although this has not yet been developed, it might provide an interesting direction for future research.
	\par 
	We expect also that the theory developed here might be useful in more algebraic contexts.
	In this direction, Chen \cite[Theorem 1.1]{ChenVolume}, cf. also Wolfe \cite[Theorem 5.14]{WolfeThes}, expressed the volumes of a relatively big line bundle on a family over a curve in terms of the limiting Harder-Narasimhan measures.
	We observe -- through the work of Kaveh-Khovanskii \cite{KavehKhov} -- that the same formula continues to hold in the non-big setting where the volume is replaced by the multiplicity. 
	Moreover, by using our version of the Mehta-Ramanathan theorem, we show that a part of the formula continues to hold even for families over higher dimensional bases.
	\par 
	Let us now proceed to the formulation of our results.
	To set up the stage, recall that a \textit{slope} (or $[\omega_B]$-slope) of a coherent sheaf $\mathscr{E}$ on $B$ is defined as 
	\begin{equation}
		\mu(\mathscr{E}) := \frac{\deg(\mathscr{E})}{\rk{\mathscr{E}}},
	\end{equation}
	where the degree, $\deg(\mathscr{E})$, is given by 
	\begin{equation}
		\deg(\mathscr{E}) := \int_B c_1(\det (\mathscr{E})) \cdot [\omega_{B, 1}] \cdots [\omega_{B, m - 1}],
	\end{equation}
	where $\det \mathscr{E}$ is Knudsen-Mumford determinant of $\mathscr{E}$, see \cite{Knudsen1976}.
	A torsion-free coherent sheaf $\mathscr{E}$ is called \textit{semistable} (or $[\omega_B]$-\textit{semistable}) if for every coherent subsheaf $\mathcal{F}$ of $\mathscr{E}$, verifying $\rk{\mathcal{F}} > 0$, we have $\mu(\mathcal{F}) \leq \mu(\mathscr{E})$.
	When the multipolarization is \textit{symmetric}, i.e. $[\omega_{B, 1}] = \cdots = [\omega_{B, m - 1}]$, we recover the usual notion of slopes and semistability, as in \cite{KobaVB}.
	Although the results of the present article are already new in the case of symmetric multipolarizations, we choose to state them in greater generality, as this appears natural in view of one of the applications \cite[Theorem 1.10]{FinHNII}.
	\par 
	Recall that any torsion-free coherent sheaf $\mathscr{E}$ on $(B, [\omega_B])$ admits a unique filtration by subsheaves $\mathscr{E}_i$, $i = 0, \ldots, h$, also called \textit{Harder-Narasimhan filtration}:
	\begin{equation}\label{eq_HN_filt}
		\mathscr{E} = \mathscr{E}_h \supset \mathscr{E}_{h - 1} \supset \cdots \supset \mathscr{E}_1 \supset \mathscr{E}_0 = \{0\},
	\end{equation}
	such that for any $1 \leq i \leq h$, the quotient sheaf $\mathscr{E}_i / \mathscr{E}_{i - 1}$ is the maximal semistable (torsion-free) subsheaf of $\mathscr{E} / \mathscr{E}_{i - 1}$, i.e. for any subsheaf of $\mathcal{F}$ of a (torsion-free) sheaf $\mathscr{E} / \mathscr{E}_{i - 1}$, we have $\mu(\mathcal{F}) \leq \mu(\mathscr{E}_i / \mathscr{E}_{i - 1})$ and $\rk{\mathcal{F}} \leq \rk{\mathscr{E}_i / \mathscr{E}_{i - 1}}$ if $\mu(\mathcal{F}) = \mu(\mathscr{E}_i / \mathscr{E}_{i - 1})$.
	For the proof, see \cite[Theorem 1.6.7]{HuybLehn} for symmetric multipolarization and \cite[Corollary 2.27]{GrebKebPetMov} for the general case.
	\par 
	We define the \textit{Harder-Narasimhan slopes}, $\mu_1, \ldots, \mu_{\rk{\mathscr{E}}}$ of $\mathscr{E}$, so that $\mu(\mathscr{E}_i / \mathscr{E}_{i - 1})$ appears among $\mu_1, \ldots, \mu_{\rk{\mathscr{E}}}$ exactly $\rk{\mathscr{E}_i / \mathscr{E}_{i - 1}}$ times, and the sequence $\mu_1, \ldots, \mu_{\rk{\mathscr{E}}}$ is non-increasing.
	We call $\mu_{\min} := \mu_{\rk{\mathscr{E}}}$ and $\mu_{\max} := \mu_{1}$, the minimal and the maximal slopes respectively. 
	Occasionally, when the relevant parameters are not clear from the context, we write $\mu_{\min}(\mathscr{E})$, $\mu_{[\omega_B], \min}(\mathscr{E})$, $\mu_i(\mathscr{E})$, etc., to denote these quantities.
	\par 
	By the flatness of $\pi$, the coherent sheaves $\mathscr{E}_k$ are torsion-free, cf. Section \ref{sect_sm_bnd}.
	We denote by $N_k := \rk{\mathscr{E}_k}$ the rank of the general fiber of $\mathscr{E}_k$, and let $\mu_1^k, \ldots, \mu_{N_k}^k$ be the Harder-Narasimhan slopes of $\mathscr{E}_k$. Let $\mu_{\min}^k, \mu_{\max}^k$ be the minimal and the maximal slopes.
	\par 
	We assume from now on that $L$ is \textit{effective on the general fiber}, i.e. for some $k \in \nat^*$, we have $N_k > 0$.
	Note that for any $k, l \in \nat$, we have a natural multiplication map 
	\begin{equation}\label{eq_sect_r_fibered_mult}
		\mathscr{E}_k \otimes \mathscr{E}_l \to \mathscr{E}_{k + l}.
	\end{equation}
	Moreover, for any $b \in B$, there is a natural map from the fiber of $\mathscr{E}_k$ at $b$ to  $H^0(X_b, L_b^k)$, where $i_b : X_b \to X$ is the fiber of $\pi$ at $b$, and $L_b := i_b^* L$.
	For general $b$, this map is an isomorphism, cf. Section \ref{sect_sm_bnd}, and the map (\ref{eq_sect_r_fibered_mult}) corresponds to the multiplication in the respective section ring.
	Since the general fiber of $\pi$ is irreducible, the section ring has no zero divisors, and so the set of $k \in \nat$ so that $N_k > 0$ forms a non-empty semigroup.
	By replacing $L$ by its tensor power associated with the greatest common divisor of this semigroup, called \textit{(relative) exponent}, we will assume without loosing the generality that $N_k > 0$ for all sufficiently big $k$.
	\par 
	Note that when $L$ is relatively ample (or even relatively big, i.e. for a certain ample line bundle $A$ on $B$, the line bundle $\pi^*A \otimes L$ is big), then the above procedure is unnecessary, as the relative exponent equals $1$, cf. \cite[Corollaries 2.2.7 and 2.2.11]{LazarBookI}.
	\par 
	For $k$ big enough, we define the probability measure $\eta_k^{HN}$ on $\real$ as follows
	\begin{equation}\label{eq_eta_defn}
		\eta_k^{HN} := \frac{1}{N_k} \sum_{i = 1}^{N_k} \delta \Big[ \frac{\mu_i^k}{k} \Big],
	\end{equation}
	where $\delta[x]$ is the Dirac mass at $x \in \real$. In Section \ref{sect_sm_bnd}, we establish the following result, serving as a starting point of the current work.
	\begin{thm}\label{thm_conv_meas}
		The measures $\eta_k^{HN}$ converge vaguely, as $k \to \infty$, to a probability measure $\eta^{HN}$ on $\real$, which has a bounded from above support, and which is absolutely continuous with respect to the Lebesgue measure, except perhaps for a point mass at $\esssup \eta^{HN}$.
		Also, the limit below exists and relates with $\eta^{HN}$ as follows 
		\begin{equation}\label{eq_conv_meas}
			\eta_{\max}^{HN} := \lim_{k \to \infty} \frac{\mu_{\max}^k}{k} = \esssup \eta^{HN}.
		\end{equation}
	\end{thm}
	\begin{rem}
	Theorem \ref{thm_conv_meas} in fact holds for a much broader class of degree functions.
		It suffices that the latter could be bounded by those associated with a multipolarization, that the notion of degree admits Harder-Narasimhan filtrations, that the analogue of Lemma \ref{lem_slope_zero_morph} holds, and that semistability behaves well under tensor products.
		Our main result, Theorem \ref{thm_mehta_ramanathan}, however, appears to be meaningful only for degrees arising from a multipolarization, which justifies the level of generality adopted in the present article.
	\end{rem}
	\par 
	When $\dim B = 1$, the convergence part of Theorem \ref{thm_conv_meas} was established by Chen \cite{ChenHNolyg} for relatively ample line bundles and then in \cite[\S 4]{ChenArithmFujita} for relatively big line bundles. 
	For this, he verified that the Harder-Narasimhan filtrations of $\mathscr{E}_k$ are submultiplicative, see (\ref{eq_sm_cond}), and that $\mu_{\max}^k$ grows at most linearly in $k$.
	Then Theorem \ref{thm_conv_meas} follows from the general theory of bounded submultiplicative filtrations developed by Chen \cite{ChenHNolyg} and Boucksom-Chen \cite{BouckChen}, cf. Theorem \ref{thm_filt}.
	\par 
	A combination of the works of Kaveh-Khovanskii \cite[Theorem 3.4]{KavehKhov} and Chen \cite{ChenArithmFujita} shows that the general theory of bounded submultiplicative filtrations extends beyond the big setting.
	Our proof of Theorem \ref{thm_conv_meas} then proceeds along the same lines as in \cite{ChenHNolyg}: by establishing the submultiplicativity of the Harder-Narasimhan filtrations of $\mathscr{E}_k$ and the linear upper bound on $\mu_{\max}^k$.
	Apart from some minor technicalities, the proof of the submultiplicativity remains the same, and the proof of the linear upper bound on $\mu_{\max}^k$ for $\dim B > 1$ reduces to the case of $\dim B = 1$.
	\par 
	Note that although Theorem \ref{thm_conv_meas} provides a thorough analysis of the asymptotic behavior of the maximal slopes, the corresponding investigation of the minimal slopes remains somewhat more delicate in our level of generality.
	In particular, while expected, it is yet unclear if 
	\begin{equation}\label{eq_add_cond}
		\liminf_{k \to \infty} \frac{\mu_{\min}^k}{k} \neq -\infty.
	\end{equation}
	It will follow from the results of Section \ref{sect_sm_bnd} that (\ref{eq_add_cond}) holds if the section ring of the restriction of $L$ to a general fiber of $\pi$ is finitely generated (in this case the liminf above is actually a limit), see Propositions \ref{prop_superadd} and \ref{prop_sm}. 
	But even in this case it is not clear if the following equality holds
	\begin{equation}\label{eq_conv_meas0}
		\eta_{\min}^{HN} := \lim_{k \to \infty} \frac{\mu_{\min}^k}{k} = \essinf \eta^{HN}.
	\end{equation}
	However, if we additionally assume that $L$ is relatively ample and $\dim B = 1$, (\ref{eq_conv_meas0}) is known to hold, see Proposition \ref{prop_xu_zh}.
	Note also that the finite generation of the section ring of the general fiber is assured, for example, when $L$ is relatively semi-ample and the general fiber of $\pi$ is normal, cf. \cite[Example 2.1.29]{LazarBookI}.
	\par 
	Already from the above discussion, we see that the measure $\eta^{HN}$ is easier to study when the base of the family is a curve.
	This is for various reasons: the subsheaves in the Harder-Narasimhan filtration become locally free; there is a very precise relation between the slope and the number of holomorphic sections, cf. Proposition \ref{prop_gl_sect_rs}, and -- finally -- the notion of slope doesn't depend on the choice of multipolarization.
	\par 
	The first main goal of this article is to show that $\eta^{HN}$, $\eta_{\max}^{HN}$ can be recovered from the similar quantities associated with base changes of our family over general curves in $B$. 
	\par 
	When $m \geq 2$ and that the multipolarization $[\omega_B]$ is \textit{integral} and \textit{very ample}, i.e. $[\omega_{B, i}] \in H^2(B, \integ)$ and very ample for any $i = 1, \ldots, m - 1$.
	Consider a regular curve $\iota: C = C(l) \hookrightarrow B$ given by Bertini theorem as an intersection of $m - 1$ general divisors in the classes $l_i  [\omega_{B, i}]$, $i = 1, \ldots, m - 1$, where $l \in \nat^{*(m - 1)}$, $l = (l_1, \ldots, l_{m - 1})$.
	It is well-known, cf. Corollary \ref{cor_welldef_slo}, that for such \textit{general} curves, the Harder-Narasimhan slopes of $\iota^* \mathscr{E}_k$ do not depend on the choice of $C$ (for a fixed parameter $l \in \nat^{*(m - 1)}$).
	Define
	\begin{equation}\label{eq_eta_cmax_edefe}
		\eta_{k, l}^{HN|C} := \frac{1}{N_k} \sum_{i = 1}^{N_k} \delta \Big[ \frac{\mu_i(\iota^* \mathscr{E}_k)}{k \cdot l_1 \cdots l_{m - 1}} \Big],
		\qquad
		\eta_{\max, k, l}^{HN|C} := \frac{\mu_{\max}(\iota^* \mathscr{E}_k)}{k \cdot l_1 \cdots l_{m - 1}}.
	\end{equation}
	For the rationale behind the normalization by $l_1 \cdots l_{m - 1}$, see (\ref{eq_slope_rest}).
	\par 
	Similarly to Theorem \ref{thm_conv_meas}, we show in Proposition \ref{prop_vague_rest_convv} that the measures  $\eta_{k, l}^{HN|C}$ converge vaguely, as $k \to \infty$, towards a probability measure $\eta_l^{HN|C}$, and the limit below exists
	\begin{equation}\label{eq_conv_measaaaasss}
		\eta_{\max, l}^{HN|C} := \lim_{k \to \infty} \eta_{\max, k, l}^{HN|C}.
	\end{equation}
	Moreover, in the setting described after (\ref{eq_conv_meas0}), the limiting minimal slopes $\eta_{\min, l}^{HN|C}$ defined analogously are also well-defined.
	\begin{thm}\label{thm_mehta_ramanathan}
		Assume that $m \geq 2$, then the measures $\eta_l^{HN|C}$ converge vaguely to $\eta^{HN}$, as $l = (l_1, \ldots, l_{m - 1}) \to \infty$ (i.e. $l_i \to \infty$, $i = 1, \ldots, m - 1$).
		Moreover,
		\begin{equation}\label{eq_mehta_ramanathan_minmax}
			\lim_{l \to \infty} \eta_{\max, l}^{HN|C} = \eta_{\max}^{HN}.
		\end{equation}
	\end{thm}
	\begin{rem}
		When $L$ is relatively ample, we moreover have
		\begin{equation}\label{eq_mehta_ramanathan_minmax0}
			\lim_{l \to \infty} \eta_{\min, l}^{HN|C} = \eta_{\min}^{HN}.
		\end{equation}
	\end{rem} 
	\par 
	As an application of our results, we show that in a fibered setting, the multiplicity of a line bundle can be estimated in terms of the Harder-Narasimhan filtration of its direct images and the multiplicity of the restriction of the line bundle to a general fiber.
	\par 
	To state our results, recall first that the volume of a line bundle $L$ over a complex projective irreducible variety $Y$ of dimension $n$ is defined as 
	\begin{equation}\label{eq_fujita_thm}
		{\rm{vol}}(L) :=
		\lim_{k \to \infty}
		\frac{\dim H^0(Y, L^{k})}{k^{n}/n!}.
	\end{equation}
	The existence of the above limit is assured by the Fujita's theorem, cf. \cite{Fujita}.
	\par 
	When ${\rm{vol}}(L) > 0$, the line bundle $L$ is called \textit{big}. 
	When $L$ is not big, the asymptotic study of $H^0(Y, L^{k})$ appears non-trivial only if $L$ is \textit{effective}, i.e. $H^0(Y, L^k) \neq \{0\}$ for some $k \in \nat^*$.
	We assume this from now on and, as described before (\ref{eq_eta_defn}), upon restricting to a sufficiently big tensor power of $L$, we further assume that $H^0(Y, L^k) \neq \{0\}$ for $k \in \nat$ big enough.
	Recall that Kodaira-Iitaka dimension, $\kappa(L) \in \{0, 1, \ldots, n\}$, of $L$ is defined, see Iitaka \cite{Iitaka}, by the property that there are $c, C > 0$, so that for $k \in \nat$ big enough, we have 
	\begin{equation}\label{eq_kod_dim_def}
		c \cdot k^{\kappa(L)}
		\leq
		\dim H^0(Y, L^{k}) 
		\leq 
		C \cdot k^{\kappa(L)}.
	\end{equation}
	While the existence of $\kappa(L)$ is an old and classical result, it was only recently established --  using the theory of Okounkov bodies -- by Kaveh-Khovanskii in \cite[Theorem 3.4]{KavehKhov} that the limit
	\begin{equation}\label{eq_exist_volume}
		{\rm{vol}}_{\kappa}(L)
		:=
		\lim_{k \to \infty} \frac{\dim H^0(Y, L^{k})}{k^{\kappa(L)} / \kappa(L)!}
	\end{equation}
	exists and is positive.
	We call it the multiplicity of $L$.
	\par 
	We now come back to our fibered setting.
	Note that as $\mathscr{E}_k$ are coherent sheaves, the Kodaira-Iitaka dimension of the general fiber, $\kappa(L_f) \in \{0, 1, \ldots, n\}$, and the multiplicity of the general fiber, ${\rm{vol}}_{\kappa}(L_f) > 0$, are well-defined, as they can be read off from the rank of $\mathscr{E}_k$.
	It then follows rather easily from the general theory of the Kodaira-Iitaka dimension, cf. Section \ref{sect_vol_est}, that
	\begin{equation}\label{eq_kappa_ineq}
		\kappa(L) \leq \kappa(L_f) + m.
	\end{equation}
	\par 
	Our last theorem provides a more precise version of (\ref{eq_kappa_ineq}) on the level of multiplicities.
	We denote by $\eta^{HN}$ the measure associated with the above family from Theorem \ref{thm_conv_meas}.
	Let us first recall that when $m = 1$ and $L$ is relatively big, Chen in \cite{ChenVolume} established that
 	\begin{equation}
	 	{\rm{vol}}(L) 
		=
		\frac{(n+1)!}{n!}  \cdot
		{\rm{vol}}(L_f)
		\cdot
		\int_{x = 0}^{+\infty} x \cdot d \eta^{HN}(x).
	\end{equation}
	By repeating his argument verbatim, while replacing the use of the existence of the limit (\ref{eq_fujita_thm}) by (\ref{eq_exist_volume}), one gets that when $m = 1$, we have
	\begin{equation}
		\lim_{k \to \infty} \frac{\dim H^0(X, L^{k})}{k^{\kappa(L_f) + m}}
		=
		\frac{{\rm{vol}}_{\kappa}(L_f)}{\kappa(L_f)!}  \cdot
		\int_{x = 0}^{+\infty} x \cdot d \eta^{HN}(x).
	\end{equation}
	\par 
	The next result shows that a part of this statement continues to hold over higher-dimensional bases.
	\begin{thm}\label{thm_vol_bnd}
		Assume that $m \geq 2$, then the following bound holds
		\begin{equation}\label{eq_vol_bnd}
			\lim_{k \to \infty} \frac{\dim H^0(X, L^{k})}{k^{\kappa(L_f) + m}}
			\leq
			\frac{{\rm{vol}}_{\kappa}(L_f)}{\kappa(L_f)! m! \cdot \prod_{i = 1}^{m - 1} a_i} 
			\cdot
			\int_{x = 0}^{+\infty} x^m \cdot d \eta^{HN}(x),
		\end{equation}
		where $a_i := \int_B [\omega_{B, 1}] \cdot \ldots \cdot [\omega_{B, m - 1}] \cdot [\omega_{B, i}]$.
		In particular, if we have $\esssup \eta^{HN} \leq 0$, then a strict inequality holds in (\ref{eq_kappa_ineq}). 
		Moreover, if $\esssup \eta^{HN} < 0$, then $L$ is not effective. 
	\end{thm}
	\begin{rem}\label{rem_mult_chen}
		We provide an explicit example in Section \ref{sect_vol_est}, showing that an equality does not necessarily hold in (\ref{eq_vol_bnd}) if $m \geq 2$.
	\end{rem}
	\par 
	Let us finish the introduction with a few words concerning our techniques.
	To prove Theorem \ref{thm_mehta_ramanathan}, we use a theorem of Mehta-Ramanathan \cite{MehtaRama}, showing that semistability behaves well under pull-backs over general curves of sufficiently big degrees. 
	Using this theorem, in (\ref{eq_thm_mehta_ramanathan_3}), we show that the coupling of $\eta^{HN}$ with an arbitrary test function on $\real$ can be calculated as a double limit, where the first limit is taken in $l$ (corresponding to the degree of the general curve, $C$) and the second limit is in $k$ (corresponding to the index of the direct image sheaves, $\mathscr{E}_k$).
	The main content of Theorem \ref{thm_mehta_ramanathan} is that we can interchange the order of these limits.
	Philosophically, even though we do not literally pursue this direction, this means that there is a certain weak uniformity in the degree of the general curve required in Mehta-Ramanathan theorem for the vector bundles $\mathscr{E}_k$ for all parameters $k \in \nat$.
	In order to show that we can change the order of limits, we rely on the techniques of Shatz \cite{Shatz}, to describe how Harder-Narasimhan slopes behave under the variation of the manifold (in our case, a specialization of the general curve) and on the analysis of Chen \cite{ChenArithmFujita} on the study of submultiplicative filtrations.
	\par 
	This article is organized as follows.
	In Section \ref{sect_gen_th_sf}, we recall some facts from the theory of submultiplicative filtrations.
	In Section \ref{sect_sm_bnd}, we establish Theorem \ref{thm_conv_meas}.
	In Section \ref{sect_deneric_slopes}, we define the general slopes of a vector bundle, which are the Harder-Narasimhan slopes of pull-backs of this vector bundle to general curves of a given degree. 
	In Section \ref{sect_deg_fam}, we study how these general slopes depend on the degree of the curve. 
	Until there, the theory we develop is for general vector bundles and hence independent of the fibration structure.
	We specialize it in Section \ref{sect_MR} to establish Theorem \ref{thm_mehta_ramanathan}.
	In Section \ref{sect_vol_est}, we establish Theorem \ref{thm_vol_bnd}.	
	In Section \ref{sect_num_char}, we establish an algebraic interpretation of $\eta_{\min}^{HN}$ and $\eta_{\max}^{HN}$.
	Finally, in Section \ref{sect_open} we discuss some possible directions for further research.
	\par
	\textbf{Notation}.
	We say that a property is satisfied for \textit{general} points in some connected complex manifold it is is satisfied for all points outside a countable union of proper analytic subspaces.
	We say that a curve $C \subset B$ is a complete intersection in an integral multipolarization $[\omega_B] = ([\omega_{B, 1}], \ldots, [\omega_{B, m - 1}])$, if it is given by the intersection of $m - 1$ divisors in $B$ in the classes $[\omega_{B, i}]$, $i = 1, \ldots, m - 1$.
	For $l = (l_1, \ldots, l_{m - 1}) \in \nat^{*(m - 1)}$, we denote $l [\omega_B] := (l_1 [\omega_{B, 1}], \ldots, l_{m - 1}[\omega_{B, m - 1}])$.
	By an abuse of notations, a $m - 1$ tuple of ample line bundles on $B$ will also be called a multipolarization.
	\par 
	A sequence $a_n \in \nat$, $n \in \nat$, is called \textit{multiplicative} if there is $b_n \in \nat^*$, such that $a_{n + 1} = a_n \cdot b_n$.
	A sequence $a_n \in \real$, $n \in \nat$, is called \textit{subadditive} (resp. \textit{superadditive}) if for any $n, m \in \nat$, we have $a_{n + m} \leq a_n + a_m$ (resp. $a_{n + m} \geq a_n + a_m$).
	We say that a property holds \textit{eventually} if it holds for sufficiently big parameters. 
	\par 
	\textbf{Acknowledgement}. 
	I would like to thank CNRS and École Polytechnique for their support and the anonymous referee for the helpful comments.
	A part of this paper was written in the Fall of 2023 during a visit in Columbia University. 
	I would like to thank the mathematical department of Columbia University, especially Duong H. Phong, for their hospitality and Alliance Program for their support.
	Since the first version of this article appeared, I have benefited from several insightful discussions with Huayi Chen and Masafumi Hattori, for which I am sincerely grateful.
	
	\section{General theory of submultiplicative filtrations}\label{sect_gen_th_sf}
	The main goal of this section is to gather some results from the theory of submultiplicative filtrations on section rings.
	\par 
	Recall that a \textit{non-increasing $\real$-filtration} $\mathcal{F}$ of a vector space $V$ is a map from $\real$ to vector subspaces of $V$, $t \mapsto \mathcal{F}^t V$, verifying $\mathcal{F}^t V \subset \mathcal{F}^s V$ for $t > s$, and such that $\mathcal{F}^t V  = V$ for sufficiently small $t$ and $\mathcal{F}^t V = \{0\}$ for sufficiently big $t$.
	We say that $\mathcal{F}$ is \textit{graded} if it respects the grading of $V$.
	It is \textit{left-continuous} if for any $t \in \real$, there is $\epsilon_0 > 0$, such that $\mathcal{F}^t V = \mathcal{F}^{t - \epsilon} V $ for any $0 \leq \epsilon \leq \epsilon_0$.
	All filtrations in this article are assumed to be non-increasing left-continuous and graded if applicable.
	We associate with any filtration the probability measure, $\eta_{\mathcal{F}}$, on $\real$, which we also call the \textit{jumping measure}, defined as follows
	\begin{equation}\label{eq_jump_meas_gen}
		\eta_{\mathcal{F}} := \frac{1}{\dim V} \sum_{j = 1}^{\dim V} \delta \big[ \eta_{\mathcal{F}}(j) \big], 
	\end{equation}
	where $\eta_{\mathcal{F}}(j)$ are the \textit{jumping numbers}, defined as $\eta_{\mathcal{F}}(j) := \sup \{ t \in \real : \dim \mathcal{F}^t V \geq j \}$.
	\par 
	Now, for a holomorphic line bundle $L$ over a complex projective irreducible variety $Y$, define the \textit{section ring}
	\begin{equation}
		R(Y, L) := \oplus_{k = 0}^{\infty} H^0(Y, L^k).
	\end{equation}
	As in (\ref{eq_kod_dim_def}), we assume from now on that $L$ is effective, and we further reduce to the case when $\dim H^0(Y, L^k) > 0$ for $k \in \nat$ big enough.
	\par 
	To transport several results from the ample to the general setting, we will rely on a result of Kaveh-Khovanskii \cite[Theorem 3.4]{KavehKhov}, who proved that the section ring $R(Y, L)$ verifies the so-called Fujita-type approximation.
	To state this result, for any $d \in \nat^*$, we consider the subring $[\sym H^0(Y, L^d)] \subset R(Y, L)$ generated by $H^0(Y, L^d)$. 
	Let $[\sym^k H^0(Y, L^d)]$ be the $k$-th graded component of $[\sym H^0(Y, L^d)]$. 
	The approximation theorem then goes as follows.
	\begin{thm}\label{thm_sym_subalgebra}
		For any $\epsilon > 0$, there is $d_0 \in \nat$ such that for any $d \geq d_0$, there is $k_0 \in \nat$ such that for any $k \geq k_0$, we have
		\begin{equation}
			\frac{\dim [\sym^k H^0(Y, L^d)]}{\dim H^0(Y, L^{kd})}
			\geq 
			1 - \epsilon.
		\end{equation}
	\end{thm}
	\par 
	We will also rely on the well-known observation, cf. \cite[Lemma 3.2]{ChangJow}, \cite[\S 3.2]{FinKIDim}, that the ring $[\sym H^0(Y, L^d)]$ coincides in high degrees with the section ring of an ample line bundle. 
	For this, recall that we have the rational maps $Y \dashrightarrow \mathbb{P}(H^0(Y, L^d)^*)$ that send a point $x \in Y$ to the hyperplane $H_{d, x} \subset H^0(Y, L^d)$ consisting of all sections in $H^0(Y, L^d)$ that vanish at $x$ (and are undefined whenever all sections in $H^0(Y, L^d)$ vanish at $x$).
	We denote by $Y_d$ the closure of the image of this map.
	It is well-known, cf. \cite[Corollary 2.1.38]{LazarBookI}, that 
	\begin{equation}\label{eq_dim_ki_yk}
		\kappa(L) = \dim(Y_d), \quad \text{for } d \in \nat \text{ large enough}.
	\end{equation}
	We denote by $A_d$ the pull-back of the hyperplane bundle $\mathscr{O}(1)$ over $\mathbb{P}(H^0(Y, L^d)^*)$ to $Y_d$.
	Then, cf. \cite[\S 3.2]{FinKIDim}, for sufficiently large $k \in \nat$, there is an isomorphism 
	\begin{equation}\label{eq_sym_yk_rel}
		[\sym^k H^0(Y, L^d)]
		\simeq
		H^0(Y_d, A_d^{k}).
	\end{equation}
	Note also that since $Y$ is irreducible, its image under the given map is also irreducible. 
	Therefore, its closure, $Y_d$, remains irreducible as well.
	\par  
	Now, a graded filtration $\mathcal{F}$ on $R(Y, L)$ is called \textit{submultiplicative} if for any $t, s \in \real$, $k, l \in \nat$, 
	\begin{equation}\label{eq_sm_cond}
		\mathcal{F}^t H^0(Y, L^k) \cdot \mathcal{F}^s H^0(Y, L^l) \subset \mathcal{F}^{t + s} H^0(Y, L^{k + l}).
	\end{equation}
	We say that $\mathcal{F}$ is \textit{bounded} if there is $C > 0$, such that for any $k \in \nat^*$, $\mathcal{F}^{ C k} H^0(Y, L^k) = \{0\}$.
	Remark that if $R(Y, L)$ is finitely generated, then there is $C > 0$, such that 
	\begin{equation}\label{eq_bnd_blow}
		\mathcal{F}^{ - C k} H^0(Y, L^k) = H^0(Y, L^k).
	\end{equation}
	For general $L$, however, this bound is not guaranteed.
	\par 
	In the context of filtrations associated with test configurations, the importance of such bounds in analytic approaches to the study of submultiplicative filtrations has been underlined by Phong-Sturm in \cite[Lemma 4]{PhongSturmTestGeodK}.
	\par 
	We let $N_k = \dim H^0(Y, L^k)$, and denote by $\eta_{\mathcal{F}}(1, k) \geq \cdots \geq \eta_{\mathcal{F}}(N_k, k)$ the \textit{jumping numbers} of the filtration $\mathcal{F}$ on $H^0(Y, L^k)$.
	We denote by $\eta_{\mathcal{F}, \max}(k) := \eta_{\mathcal{F}}(1, k)$ the maximal and by $\eta_{\mathcal{F}, \min}(k) := \eta_{\mathcal{F}}(N_k, k)$ the minimal jumping numbers.
	The following result is well-known in the ample case, cf. \cite[Propositions 3.2.4, 3.2.6]{ChenHNolyg}.
	The proof in the general case remains verbatim.
	\begin{prop}\label{prop_superadd}
		For any submultiplicative filtration $\mathcal{F}$ on $R(Y, L)$, the sequence $\eta_{\mathcal{F}, \max}(k)$, $k \in \nat$, is superadditive.
		In particular, the following limit exists 
		\begin{equation}
			\eta_{\mathcal{F}, \max} := \lim_{k \to \infty} \frac{\eta_{\mathcal{F}, \max}(k)}{k}, 
		\end{equation}
		in $\real \cup \{ + \infty \}$, and it is finite if $\mathcal{F}$ is bounded.
		If $R(Y, L)$ is finitely generated, the sequence $\eta_{\mathcal{F}, \min}(k)$ is eventually superadditive, and the limit $\lim_{k \to \infty} \eta_{\mathcal{F}, \min}(k) / k$ exists in $\real \cup \{ + \infty \}$, and it is finite if $\mathcal{F}$ is bounded.
	\end{prop}
	Define the sequence of \textit{jumping measures} $\eta_{\mathcal{F}, k}$, $k \in \nat^*$, of $\mathcal{F}$ on $\real$ as follows
	\begin{equation}\label{eq_jump_meas_d}
		\eta_{\mathcal{F}, k} := \frac{1}{N_k} \sum_{j = 1}^{N_k} \delta \Big[ \frac{\eta_{\mathcal{F}}(j, k)}{k} \Big].
	\end{equation}
	\par 
	For any $d \in \nat^*$, we define the filtration $\mathcal{F}^d$ on $[{\rm{Sym}} H^0(Y, L^d)]$ as the maximal submultiplicative filtration which coincides with the restriction of the filtration $\mathcal{F}$ on $H^0(Y, L^d)$.
	As $[{\rm{Sym}} H^0(Y, L^d)]$ is generated by $H^0(Y, L^d)$, the resulting filtration is clearly bounded.
	We denote by $\eta_{\mathcal{F}, k}^{d}$ the associated probability measure on $\real$, defined as 
	\begin{equation}
		\eta_{\mathcal{F}, k}^{d}
		:=
		\frac{1}{\dim [{\rm{Sym}}^k H^0(Y, L^d)]} \sum_{j = 1}^{\dim [{\rm{Sym}}^k H^0(Y, L^d)]} \delta \big[ \frac{\eta_{\mathcal{F}^d}(j, k)}{k d} \big],
	\end{equation}
	\par 
	For any $\lambda \in \real$, consider the subspace $R^{(\lambda)}(Y, L) := \oplus_{k = 0}^{\infty} \mathcal{F}^{\lambda k} H^0(Y, L^k)$ of $R(Y, L)$.
	Note that it forms a subalgebra by the submultiplicativity assumption on $\mathcal{F}$.
	\par 
	\begin{thm}\label{thm_filt}
		For any bounded submultiplicative filtration $\mathcal{F}$ on $R(Y, L)$, the jumping measures $\eta_{\mathcal{F}, k}$, $k \in \nat^*$, converge vaguely, as $k \to \infty$, to a probability measure $\eta_{\mathcal{F}}$ on $\real$ with a bounded from above support.
		Moreover, for any $d \in \nat^*$, $\eta_{\mathcal{F}, k}^{d}$ converge weakly, as $k \to \infty$, to a probability measure $\eta_{\mathcal{F}}^d$ on $\real$ with a compact support, and as $d \to \infty$, the sequence of measures $\eta_{\mathcal{F}}^{d}$ converge vaguely towards $\eta_{\mathcal{F}}$.
		Also, the following identity holds
		\begin{equation}\label{eq_filt}
			\lim_{k \to \infty} \frac{\eta_{\mathcal{F}, \max}(k)}{k} = \esssup \eta_{\mathcal{F}},
		\end{equation}
		and for any $\lambda \in \real$, the following limit exists
		\begin{equation}\label{eq_vol_conc}
			{\rm{vol}} (R^{(\lambda)}(Y, L)) := \lim_{k \to \infty} \frac{ \dim \mathcal{F}^{\lambda k} H^0(Y, L^k)}{ N_k } = \int_{\lambda}^{+ \infty} d \eta_{\mathcal{F}}.
		\end{equation}
		Finally, if $\kappa(L) > 0$, then $\lambda \mapsto {\rm{vol}} (R^{(\lambda)}(Y, L))^{1 / \kappa(L)}$ is concave on $]-\infty, \esssup \eta_{\mathcal{F}}[$. 
	\end{thm}
	\begin{rem}
		a)
		Boucksom-Chen in \cite{BouckChen} gave an explicit expression of $\eta_{\mathcal{F}}$ in terms of the Okounkov body associated with $\mathcal{F}$.
		Witt Nyström \cite{NystOkounTest} gave an expression of $\eta_{\mathcal{F}}$ in terms of the geodesic ray in the space of positive metrics on $L$ associated with the filtration, when the filtration is given by the weight of a $\comp^*$-action; Hisamoto \cite{HisamSpecMeas} later generalized it for finitely generated filtrations; and the author in \cite{FinNarSim} further generalized it to all bounded submultiplicative filtrations by giving a different proof, based purely on complex pluripotential theory and semiclassical input from Tian \cite{TianBerg} and Phong-Sturm \cite{PhongSturm}.
		This was further generalized by the author in \cite{FinTits} in the relative setting, following Chen-Maclean \cite{ChenMaclean}, cf. also Boucksom-Jonsson \cite{BouckJohn21}, in the big setting in \cite{FinGQBig} and in the pointwise sense, i.e. for the respective weighted Bergman kernels in \cite{FinSubmToepl}.
		\par 
		b) While we obviously have the inequality $\limsup_{k \to \infty} \frac{\eta_{\mathcal{F}, \min}(k)}{k} \leq \essinf \eta_{\mathcal{F}}$, the equality does not always hold, cf. \cite[Remark 5.5]{BouckHisJonsFourier}.
		\par 
		c) Note that the last statement implies that the measure $\eta_{\mathcal{F}}$ is absolutely continuous with respect to the Lebesgue measure, except perhaps for a point mass at $\esssup \eta_{\mathcal{F}}$.
		\par 
		d) The statement and the argument below remain valid even if the filtration $\mathcal{F}$ is defined only over a \textit{linear series} (i.e. a graded subring) $W = \oplus_{k = 0}^{\infty} W_k \subset R(Y, L)$, see the discussion around \cite[(3.16)]{FinKIDim} for a related analysis.
	\end{rem}
	\begin{proof}
		When the line bundle is ample, the weak convergence of jumping measures was established by Chen \cite{ChenHNolyg} and Boucksom-Chen \cite{BouckChen}, cf. \cite[Theorem 5.3]{BouckHisJonsFourier}.
		They also verified that (\ref{eq_filt}) and (\ref{eq_vol_conc}) hold in that setting.
		In \cite[Theorem 2.11]{ChenArithmFujita}, Chen also established the vague convergence of jumping measures on arbitrary algebras, which are \textit{approximable}.
		Theorem \ref{thm_sym_subalgebra} says precisely that $R(Y, L)$ is approximable.
		This shows that the jumping measures $\eta_{\mathcal{F}, k}$ converge vaguely, as $k \to \infty$, to a probability measure $\eta_{\mathcal{F}}$ on $\real$.
		\par 
		By the identification (\ref{eq_sym_yk_rel}) and the fact that $Y_d$ is irreducible, the ample case of Theorem \ref{thm_filt} shows that $\eta_{\mathcal{F}, k}^d$ converge weakly towards a measure of compact support on $\real$, which we denote by $\eta_{\mathcal{F}}^d$.
		Boucksom-Jonsson in \cite[Theorems 3.12 and 3.18]{BouckJohn21} proved that the measures $\eta_{\mathcal{F}}^d$ converge vaguely, as $d \to \infty$, towards $\eta_{\mathcal{F}}$.
		Note that the authors there assumed that $L$ is ample, but this assumption was never used, except for \cite[Lemma 3.21]{BouckJohn21}, the analogue of which continues to holds without ampleness by Kaveh-Khovanskii \cite[Theorem 3.4]{KavehKhov}.
		It is hence now only left to establish (\ref{eq_filt}) and (\ref{eq_vol_conc}).
		We show that both of these statement follow from their ample analogues.
		\par 
		Immediately from the fact that $\esssup \eta_{\mathcal{F}} \geq \esssup \eta_{\mathcal{F}}^d$ and the vague convergence of $\eta_{\mathcal{F}}^d$ towards $\eta_{\mathcal{F}}$, we conclude that $\lim_{d \to \infty} \esssup \eta_{\mathcal{F}}^d = \esssup \eta_{\mathcal{F}}$.
		It follows then that (\ref{eq_filt}) holds in full generality, as it holds for $\eta_{\mathcal{F}}^d$ in place of $\eta_{\mathcal{F}}$.
		\par 
		We now assume that $d \in \nat^*$ is big enough so that (\ref{eq_dim_ki_yk}) holds.
		Note also that the function $g_d(\lambda) := (\int_{\lambda}^{+\infty} d \eta_{\mathcal{F}}^d)^{1/\kappa(L)}$ is concave on $]-\infty, \esssup \eta_{\mathcal{F}}^d[$ by the already established ample case of Theorem \ref{thm_filt}.
		Moreover, by the Portemanteau theorem, outside of a countable subset of $]-\infty, \esssup \eta_{\mathcal{F}}[$, $g_d$ converge pointwise, as $d \to \infty$, towards $g(\lambda) := (\int_{\lambda}^{+\infty} d \eta_{\mathcal{F}})^{1/\kappa(L)}$.
		It follows that $g$ is concave on $]-\infty, \esssup \eta_{\mathcal{F}}[$, as it is a limit (outside of a countable subset) of decreasing concave functions, finishing the proof.
	\end{proof}
	
	\section{Submultiplicative nature of the Harder-Narasimhan filtrations}\label{sect_sm_bnd}
	The main goal of this section is to prove Theorem \ref{thm_conv_meas} by showing that the theory from Section \ref{sect_gen_th_sf} applies to the Harder-Narasimhan filtrations.
	We use the notations from Introduction.
	\par 
	Let us first briefly recall why the $\mathscr{E}_k$ are torsion-free. 
	For this, we consider an open neighbourhood $U \subset B$ and a holomorphic function $f : U \to \comp$.
	As $B$ is regular (and hence, locally irreducible), the multiplication by $f$ yields the exact sequence of sheaves $0 \to \mathscr{O}_U \overset{f}{\to} \mathscr{O}_U$.
	We let $V := \pi^{-1}(U)$. 
	As $\pi$ is flat, we have the exact sequence of sheaves $0 \to \mathscr{O}_V \overset{\pi^* f}{\to} \mathscr{O}_V$.
	Multiplying this sequence by $L^{k}$ and forming the long exact sequence in cohomology, we obtain the exact sequence of sheaves $0 \to \mathscr{E}_k \overset{f}{\to} \mathscr{E}_k$. 
	As $f$ is arbitrary, this shows precisely that $\mathscr{E}_k$ is torsion-free.
	\par 
	Now, for a given point $b \in B$, we denote by $\kappa(b) \simeq \comp$ the residue field at $b \in B$, defined as $\kappa(b) = \mathscr{O}_{B, b} / \mathfrak{m}_b$, where $\mathscr{O}_{B, b}$ is the local ring of holomorphic functions at $b$, and $\mathfrak{m}_b \subset \mathscr{O}_{B, b}$ is the maximal ideal.
	For any $b \in B$, we let $\mathscr{E}_{k, b} \otimes_{\mathscr{O}_{B, b}} \kappa(b)$ be the fiber of $\mathscr{E}_{k}$ at $b$.
	Clearly, if we denote by $i_b : X_b \to X$ the fiber of $\pi$ at $b$, and $L_b = i_b^* L$, we have a natural map
	\begin{equation}\label{eq_map_restriction_sheaf}
		\mathscr{E}_{k, b} \otimes_{\mathscr{O}_{B, b}} \kappa(b)
		\to
		H^0(X_b, L_b^k).
	\end{equation}
	\par 
	As it is well-known, there is a proper analytic subset $\Sigma_k \subset B$, so that $\mathscr{E}_{k}$ is locally free outside of $\Sigma_k$, and the map (\ref{eq_map_restriction_sheaf}) is an isomorphism for any $b \notin \Sigma_k$, cf. \cite[Theorem III.2.18 and Corollary III.12.9]{HartsBook}.
	For further purposes, we let $\Sigma := \cup_{k = 0}^{+\infty} \Sigma_k$.
	\par 
	Note that (\ref{eq_map_restriction_sheaf}) might not be an isomorphism even at the locus where $\mathscr{E}_{k}$ is locally free, and $\Sigma_k$ can have codimension $1$.
	However, if $L$ is relatively ample, then $\Sigma_k = \emptyset$ for $k$ big enough by \cite[Corollary III.12.9]{HartsBook}, as Serre vanishing and the flatness of $\pi$ implies that the dimension of the right-hand side of (\ref{eq_map_restriction_sheaf}) is locally constant, cf. \cite[Theorem III.9.9]{HartsBook}.
	\par  
	We introduce the (non-increasing) filtrations $\mathcal{F}_k(\lambda)$, $\lambda \in \real$, of $\mathscr{E}_k$ by coherent (torsion-free) subsheaves (defined over $B$) so that $\mathcal{F}_k(\lambda)$ is the maximal subsheaf of $\mathscr{E}_k$ such that all of its Harder-Narasimhan slopes are not smaller than $\lambda$.
	In other words, we define
	\begin{equation}\label{eq_fklam_defnff}
		\mathcal{F}_k(\lambda) := \max \Big\{
			\mathcal{E} \subset \mathscr{E}_k : \text{ for any quotient sheaf } \mathcal{Q} \text{ of } \mathcal{E}, \mu(\mathcal{Q}) \geq \lambda
		\Big\}.
	\end{equation}
	Clearly, the filtration $\mathcal{F}_k$ is just a “renaming" of the Harder-Narasimhan filtration of $\mathscr{E}_k$.
	\par 
	\begin{sloppypar}
	The main idea behind the above “renaming" of the Harder-Narasimhan filtration is that the filtration then becomes submultiplicative.
	More precisely, for $b \in B \setminus \Sigma$, we denote by $\mathcal{F}_b$ the filtration induced by $\mathcal{F}_k(\lambda)$ on $R(X_b, L_b) = \oplus_{k = 0}^{\infty} H^0(X_b, L_b^{k})$ using (\ref{eq_map_restriction_sheaf}).
	\begin{prop}\label{prop_sm}
		There is a subset $\Sigma' \subset B \setminus \Sigma$, given by a countable union of proper analytic subsets so that the filtration $\mathcal{F}_b$ is submultiplicative for any $b \in B \setminus (\Sigma \cup \Sigma')$.
	\end{prop}
	When $\dim B = 1$, Proposition \ref{prop_sm} is due to Chen \cite[\S 4.3.2]{ChenHNolyg}.
	In this setting, if one additionally assumes that $L$ is relatively ample, the statement holds even for every $b \in B$.
	The proof of the general case remains verbatim modulo some technical modifications which enter into play due to the fact that Harder-Narasimhan filtrations are filtrations by subsheaves (and not by vector bundles) if $\dim B > 1$.
	More precisely, we need the following basic results.
	\end{sloppypar}
	\begin{lem}\label{lem_max_slope_HN}
		For any $\lambda \in \real$, we have $\mu_{\max}(\mathscr{E}_k / \mathcal{F}_k(\lambda) ) < \lambda$.
	\end{lem}
	\begin{proof}
		It follows directly from the definition of $\mathcal{F}_k(\lambda)$.
	\end{proof}
	For the following result, we recall that the determinant line of a coherent sheaf is canonically trivial if the support of the sheaf has codimension at least $2$, cf. \cite[Proposition 5.6.14]{KobaVB}.
	In particular, Harder-Narasimhan filtrations can be defined for coherent sheaves, which are torsion-free in codimension $2$, cf. \cite[Theorem 1.6.7]{HuybLehn}.
	In this case, the filtration will be given by coherent subsheaves, which are torsion-free in codimension $2$, and their subsequent quotients will also be torsion-free in codimension $2$.
	In particular, the notions of maximal and minimal slopes are well-defined. 
	We borrow for those notions the previously-defined notations of the maximal and the minimal Harder-Narasimhan slopes.
	\begin{lem}[{ \cite[Proposition 1.2.7  and Theorem 1.6.6]{HuybLehn} }]\label{lem_slope_zero_morph}
		Let $\phi : \mathscr{E} \to \mathcal{F}$ be a map between coherent sheaves $\mathscr{E}, \mathcal{F}$, such that both $\mathscr{E}$ and $\mathcal{F}$ are torsion-free in codimension $2$.
		Assume that $\mu_{\max}(\mathcal{F}) < \mu_{\min}(\mathscr{E})$. 
		Then $\phi = 0$ outside an analytic subset of codimension at least $2$.
	\end{lem}
	\begin{rem}
		In \cite{HuybLehn}, authors treat symmetric multipolarization, but the generalization is direct.
	\end{rem}
	The reason why it is necessary to consider sheaves which are not torsion-free in the context of this article is that a tensor product of two torsion-free coherent sheaves is not torsion-free in general, cf. \cite[p. 261]{RuppSera}. However, the following result is true.
	\begin{lem}\label{lem_tens_prod}
		Let $\mathscr{E}_1, \mathscr{E}_2$ be two torsion-free coherent sheaves, then their tensor product is torsion-free in codimension $2$. Similarly, for a torsion-free coherent sheaf $\mathscr{E}$, the sheaves $\sym^k \mathscr{E}$, $k \in \nat^*$, are torsion-free in codimension $2$.
	\end{lem}
	\begin{proof}
		The statement follows directly from the fact that torsion-free coherent sheaves are locally free in codimension $2$, cf. \cite[Theorem 5.5.8, Corollary 5.5.15]{KobaVB}
	\end{proof}
	\begin{proof}[Proof of Proposition \ref{prop_sm}]
		It suffices to verify that the natural map $\mathcal{F}_k(\lambda) \otimes \mathcal{F}_l(\mu) \to \mathscr{E}_{k + l} / \mathcal{F}_{k + l}(\lambda + \mu)$ is zero outside a subset of codimension at least $2$, for any $\lambda, \mu \in \real$, $k, l \in \nat$, as one could then take $\Sigma'$ to be the union of the supports of these maps.
		But this is a direct consequence of Lemmas \ref{lem_max_slope_HN}, \ref{lem_slope_zero_morph}, \ref{lem_tens_prod}, the fact that a tensor product of two semistable sheaves is semistable (in the category of sheaves which are torsion-free in codimension $2$), and the slope behaves additively with respect to tensor products.
		The last two results were established by Maruyama \cite{MaruyTensProd} in the case of symmetric mutipolarization and by Greb-Kebekus-Peternell \cite[Theorem 4.2]{GrebKebPetMov} in general.
	\end{proof}
	\begin{prop}\label{thm_bound}
		The filtration $\mathcal{F}_b$ is bounded for any $b \in B \setminus \Sigma$.
	\end{prop}
	When $m := \dim B = 1$, Proposition \ref{thm_bound} was established by Chen \cite[Proposition 4.3.5]{ChenHNolyg}, cf. also \cite[Lemma-Definition 2.26]{XuZhuPosCM}.
	We will show below that the statement for projective $B$ of any dimension reduces to the one considered by Chen.
	For this, let us establish the following result.
	\begin{lem}\label{lem_major_hn_restr}
		Assume that the multipolarization $[\omega_B]$ is integral and very ample.
		Let $\mathscr{E}$ be a torsion-free coherent sheaf over $B$.
		Then for general regular curves $\iota: C \hookrightarrow B$, which are complete intersections in $l [\omega_B]$, $l \in \nat^{*(m - 1)}$, the pull-back, $\iota^* \mathscr{E}$, of $\mathscr{E}$ to $C$ is torsion-free, and we have
		\begin{equation}
		\begin{aligned}
			l_1 \cdots l_{m - 1} \cdot \mu_{[\omega_B], \max}(\mathscr{E}) \leq \mu_{\max}(\iota^* \mathscr{E}),
			\\
			l_1 \cdots l_{m - 1} \cdot \mu_{[\omega_B], \min}(\mathscr{E}) \geq \mu_{\min}(\iota^* \mathscr{E}).
		\end{aligned}
		\end{equation}
	\end{lem}
	\begin{proof}
		First, the statement about the torsion-free pull-back is well-known, cf. \cite[Théorème 12.1.1]{EGA4} or \cite[Lemma 3.1.1]{HuybLehn}.
		Now, as $C$ is a complete intersection in $l [\omega_B]$, we have
		\begin{equation}\label{eq_slope_rest}
			\mu(\iota^* \mathscr{E}) = l_1 \cdots l_{m - 1} \cdot \mu_{[\omega_B]}(\mathscr{E}).
		\end{equation}
		Remark also that the maximal and the minimal slopes admit the following useful characterization
		\begin{equation}\label{eq_max_min_slope}
			\begin{aligned}
				&
				\mu_{[\omega_B], \max}(\mathscr{E})
				=
				\max \Big\{ 
					\mu_{[\omega_B]}(\mathcal{F})
					:
					\mathcal{F}
					\text{ is a subsheaf of } \mathscr{E}
				\Big\}
				,
				\\
				&
				\mu_{[\omega_B], \min}(\mathscr{E})
				=
				\min \Big\{ 
					\mu_{[\omega_B]}(\mathcal{Q})
					:
					\mathcal{Q}
					\text{ is a quotient sheaf of } \mathscr{E}
				\Big\}.
			\end{aligned}
		\end{equation}
		The result now follows from (\ref{eq_slope_rest}) and (\ref{eq_max_min_slope}) applied for $\mathscr{E}$ and $\iota^* \mathscr{E}$.
	\end{proof}	
	\par 
	Now, for a given curve $\iota: C \hookrightarrow B$, we denote by $\pi^C : X^C \to C$ (resp. $L_C$) the pull-back of $\pi$ (resp. $L$) through $\iota$.
	Remark that $\pi^C$ is flat, cf. \cite[Proposition III.9.3]{HartsBook}, and so $R^0 \pi^C_* L_C^{\otimes k}$ is locally free.
	Indeed, by the discussion in the beginning of the section, $R^0 \pi^C_* L_C^{\otimes k}$ is torsion-free, which implies that it is locally free, as $C$ is a curve.
	\par 
	Similarly to (\ref{eq_map_restriction_sheaf}), we have a natural map
	\begin{equation}\label{eq_map_rel_restr_curv_pushforward}
		\iota^* \mathscr{E}_{k}
		\to 
		R^0 \pi^C_* L_C^{\otimes k}.
	\end{equation}
	Note that the map (\ref{eq_map_restriction_sheaf}) factorizes through the analogous map induced  on the fibers from (\ref{eq_map_rel_restr_curv_pushforward}).
	So if $c \in C$ is such that $\iota(c) \notin \Sigma_k$, then the map induced by (\ref{eq_map_rel_restr_curv_pushforward}) on the fibers is an isomorphism.
	\par 
	As a general regular curve $\iota: C \hookrightarrow B$ which is a complete intersections in $l [\omega_B]$, $l \in \nat^{*(m - 1)}$ intersects $\Sigma_k$ only at a finite number of points, we obtain that there is a short exact sequence 
	\begin{equation}\label{eq_torsion_quotient_direct_images}
		0
		\to
		\iota^* \mathscr{E}_{k}
		\to 
		R^0 \pi^C_* L_C^{\otimes k}
		\to
		Q
		\to
		0,
	\end{equation}
	where $Q$ is a torsion sheaf on $C$. 
	Indeed, both the kernel and the cokernel of (\ref{eq_map_rel_restr_curv_pushforward}) are torsion sheaves by the above argument, but as $\iota^* \mathscr{E}_{k}$ is torsion-free by Lemma \ref{lem_major_hn_restr}, the kernel is empty.
	\par 
	We assume now that $C$ is such that $\iota^* \mathscr{E}_{k}$ is torsion-free for any $k \in \nat$, and $\iota(C) \not\subset \Sigma$.
	Similarly to (\ref{eq_fklam_defnff}) we introduce the (non-increasing) filtrations $\mathcal{F}_k^C(\lambda)$, $\lambda \in \real$, of $\iota^* \mathscr{E}_k$ by vector bundles (defined over $C$) so that
	\begin{equation}
		\mathcal{F}_k^C(\lambda) := \max \Big\{
			\mathcal{E} \subset \iota^* \mathscr{E}_k : \text{ for any quotient vector bundle } \mathcal{Q} \text{ of } \mathcal{E}, \mu(\mathcal{Q}) \geq \lambda
		\Big\}.
	\end{equation}
	Let $c \in C$ be such that $\iota(c) \not\in \Sigma$.
	Let $\mathcal{F}_c^C$ be the filtration induced by $\mathcal{F}_k^C(\lambda)$ on $R(X_{\iota(c)}, L_{\iota(c)})$ using (\ref{eq_map_restriction_sheaf}).
	By repeating the proof of Proposition \ref{prop_sm}, we obtain the following result.
	\begin{prop}\label{prop_sm2}
		The filtration $\mathcal{F}_c^C$ is submultiplicative.
	\end{prop}
	We will also need the following analogue of Proposition \ref{thm_bound}.
	\begin{prop}\label{thm_bound2}
		The filtration $\mathcal{F}_c^C$ is bounded.
	\end{prop}
	\begin{proof}[Proof of Propositions \ref{thm_bound} and \ref{thm_bound2}]
		When $m = 1$, Proposition \ref{thm_bound} was established by Chen \cite[Proposition 4.3.5]{ChenHNolyg}.
		We will show that both Propositions \ref{thm_bound} and \ref{thm_bound2} follow from this.
		Indeed, by (\ref{eq_torsion_quotient_direct_images}), any subsheaf of $\iota^* \mathscr{E}_k$ is also a subsheaf of $R^0 \pi^C_* L_C^{\otimes k}$.
		As the maximal slopes of $R^0 \pi^C_* L_C^{\otimes k}$ grow at most linearly by the work of Chen, the same holds for the maximal slopes of $\iota^* \mathscr{E}_k$, establishing Proposition \ref{thm_bound2}.
		\par 
		Now, since $B$ is projective, any multipolarization can be bounded from above by an integral multipolarization. Remark that the slopes will only increase by this procedure, so for the sake of proving Proposition \ref{thm_bound}, it is sufficient to assume that the multipolarization $[\omega_B]$ is integral and very ample.
		By Lemma \ref{lem_major_hn_restr}, the maximal Harder-Narasimhan $[\omega_B]$-slope of $\mathscr{E}_k$ is bounded from above by the maximal Harder-Narasimhan slope of $\iota^* \mathscr{E}_k$.
		But the latter grow at most linearly by the already established Proposition \ref{thm_bound2}, finishing the proof.
	\end{proof}
	\begin{proof}[Proof of Theorem \ref{thm_conv_meas}]
		It follows directly from Theorem \ref{thm_filt} and Propositions \ref{prop_sm}, \ref{thm_bound}.
	\end{proof}
	\begin{rem}
		In \cite[Theorem 1.5]{FinHNII}, the author gave a differential-geometric proof of Proposition \ref{thm_bound}, relying on Ma-Zhang \cite{MaZhangSuperconn}, which works more generally for Kähler $B$ (but asks additionally that $\pi$ is a submersion and that $L$ is relatively ample).
	\end{rem}
	\par 
	Proceeding as above, in the notations from Introduction, we obtain the following result.
	\begin{prop}\label{prop_vague_rest_convv}
		The measures  $\eta_{k, l}^{HN|C}$ converge vaguely, as $k \to \infty$, towards a probability measure $\eta_l^{HN|C}$.
		Also, the limit below exists and relates with $\eta_l^{HN|C}$ as follows 
		\begin{equation}\label{eq_conv_measaaa}
			\eta_{\max, l}^{HN|C} := \lim_{k \to \infty} \frac{\mu_{\max}(\iota^* \mathscr{E}_k)}{k \cdot l_1 \cdots l_{m - 1}} = \esssup \eta_l^{HN|C}.
		\end{equation}
	\end{prop}

	\section{Slopes of pull-backs to general curves}\label{sect_deneric_slopes}
	The main goal of this section is to define the Harder-Narasimhan slopes of pull-backs of a torsion-free coherent sheaf to general curves of a given degree. 
	For this, we first recall how the Harder-Narasimhan slopes vary in flat families.
	\par 
	More precisely, our setting will be as follows. Let $p : B \to S$ be a proper surjective flat holomorphic map of relative dimension $m$ between complex connected manifolds $B$ and $S$, and let $F = (F_1, \ldots, F_{m - 1})$ be a relative multipolarization, which means that $F_i$ are relatively ample line bundles with respect to $p$ on $B$, $i = 1, \ldots, m - 1$.
	We assume for simplicity that $p$ has a smooth fiber, which implies that the set of smooth fibers is an open dense subset of $S$, see \cite[Théorème 12.2.4]{EGA4}.
	Let $\mathscr{E}$ be a coherent torsion-free sheaf on $B$, which is flat over $S$, cf. \cite[Definition 2.1.1]{HuybLehn}. 
	The set of fibers of $p$ such that the pull-backs of $\mathscr{E}$ to them is torsion-free is open in $S$ by \cite[Théorème 12.2.1]{EGA4}.
	We also assume it is non-empty.
	Then for general $s \in S$, the Harder-Narasimhan filtrations (with respect to the multipolarization $i_s^* c_1(F) := (i_s^* c_1(F_1), \ldots, i_s^*  c_1(F_{m - 1}))$) of the torsion-free coherent sheaf $i_s^* \mathscr{E}$ are well-defined, where $i_s : B_s \to B$ is the fiber of $p$ at $s$.
	\begin{prop}\label{prop_cnst_generic}
		There is a proper analytic subset $A \subset S$, so that the $i_s^* c_1(F)$-slopes of the Harder-Narasimhan filtrations of $i_s^* \mathscr{E}$ are constant for $s \in S \setminus A$.
	\end{prop}
	The following two results will make the proof of Proposition \ref{prop_cnst_generic}.
	\begin{lem}[ { \cite[Example 20.3.3]{FultonIntTheory} } ]\label{lem_cnst_slope} 
		The $i_s^* c_1(F)$-slopes of $i_s^* \mathscr{E}$ are locally-constant functions of $s \in S$.
	\end{lem}
	\begin{thm}\label{thm_rel_hn}
		There is a complex manifold $T$ with a birational morphism $f : T \to S$, such that for $p_0 : B \times_{S} T \to B$, there is a filtration of $p_0^* \mathscr{E}$ by subsheaves:
		\begin{equation}
			p_0^* \mathscr{E} = HN(\mathscr{E})_h \supset HN(\mathscr{E})_{h - 1} \supset \cdots \supset HN(\mathscr{E})_1 \supset HN(\mathscr{E})_0 = \{0\},
		\end{equation}
		such that the factors $HN(\mathscr{E})_i / HN(\mathscr{E})_{i - 1}$ are flat for $i = 1, \ldots, h$, and there is a proper analytic subset $A \subset T$, such that the pull-back of the above filtration to a fiber at $t \in T \setminus A$ corresponds to the Harder-Narasimhan filtration of the fiber.
	\end{thm}
	\begin{rem}
		The filtration above is called \textit{relative Harder-Narasimhan filtration}. The proof of the above statement for symmetric multipolarizations can be found in \cite[Theorem 2.3.2]{HuybLehn}, and in the general case in \cite[Corollary 6.6]{TomaHNDegree}.
	\end{rem}
	\begin{proof}[Proof of Proposition \ref{prop_cnst_generic}]
		It follows directly from Theorem \ref{thm_rel_hn} and Lemma \ref{lem_cnst_slope}, applied to all subsheaves of the relative Harder-Narasimhan filtration.
	\end{proof}
	\par 
	We now apply the above general theory to the study Harder-Narasimhan filtrations of the pull-backs of a given sheaf to general curves. 
	Fix a complex projective manifold $B$ of dimension $m$ with very ample line bundles $F_1, \ldots, F_{m - 1}$.
	For $l = (l_1, \ldots, l_{m - 1}) \in \nat^{*(m - 1)}$, we consider
	\begin{equation}\label{eq_defn_sl_var}
		S_l := \mathbb{P}(H^0(B, F_1^{l_1})) \times \cdots \times \mathbb{P}(H^0(B, F_{m - 1}^{l_{m - 1}})),
	\end{equation}
	parametrizing $\dim B - 1$ hyperplanes in $H^0(B, F_1^{l_1})^*, \ldots, H^0(B, F_{m - 1}^{l_{m - 1}})^*$.
	We define the \textit{incidence variety} $\mathcal{C}_l$ as follows 
	\begin{equation}\label{eq_family_cl}
		\mathcal{C}_l := \Big\{ (b \in B, s \in S_l) : b \text{ lies in the intersection of hyperplanes parametrized by } s \Big\},
	\end{equation}
	where we implicitly identified $B$ with its images in $\mathbb{P}(H^0(B, F_i^{l_i})^*)$, $i = 1, \ldots, m - 1$, through the Kodaira map.
	We then have the natural maps $p_l : \mathcal{C}_l \to S_l$ and $p_{0, l} : \mathcal{C}_l \to B$.
	It is classical, cf. \cite[\S 3.1]{HuybLehn}, that $p_{0, l}$ is a locally trivial bundle with fibers given by the product of projective spaces. 
	In particular, $\mathcal{C}_l$ is a smooth manifold. 
	The fibers of the projection $p_l$ are given by the intersection of the $\dim B - 1$ divisors given by the zero-locus of sections parametrized by a point in $S_l$.
	Hence, these fibers are generically curves. 
	More precisely, the following result holds true.
	\begin{prop}\label{prop_bertini}
		The set $U_l \subset S_l$ (resp. $V_l \subset U_l$) such that the fibers of $p_l$ at $U_l$ have dimension $1$ (resp. of regular values of $p_l$) is a non-empty open subset in Zariski topology, i.e. its complement is an analytic subset.
	\end{prop}
	\begin{proof}
		The statement for $V_l$ is the Bertini theorem, cf. \cite[Theorem II.8.18]{HartsBook}.
		Let us now establish the statement about $U_l$.
		Define the function $n : S_l \to \nat$ as 
		\begin{equation}\label{eq_dim_max_irr_comp}
			n(s) = \text{maximal dimension of an irreducible component of } p_l^{-1}(s).
		\end{equation}
		Since $p_l^{-1}(s)$ is given by the intersection of $\dim B - 1$ divisors, by Serre's inequality on height, cf. \cite[Exercise II.11.11]{DemCompl}, the dimension of every irreducible component of $p_l^{-1}(s)$, $s \in S_l$, is at least $1$.
		In particular, we have $U_l := \{ s \in S_l : n(s) < 2 \}$.
		But $\mathcal{C}_l$ is projective and, hence, $p_l$ is proper.
		In particular, the function $n$ is upper semicontinuous (in Zariski topology), cf. \cite[Théorème 13.1.5]{EGA4}. 
		Hence, the set $\{ s \in S_l : n(s) < 2 \}$ is open, which finishes the proof. 
	\end{proof}
	\begin{prop}\label{prop_flat_fam_inc}
		The restriction of $p_l$ to $p_l^{-1}(U_l)$ is flat.
	\end{prop}
	\begin{rem}
		The same statement appeared in \cite[Proposition 1.5]{MehtaRamMathAnn} with a different proof.
	\end{rem}
	\begin{proof}
		Both $U_l$ and $p_l^{-1}(U_l)$ are smooth manifolds.
		The dimension of the fibers of $p_l$ over $U_l$ equals to $1$ by the definition of the set $U_l$, and $\dim \mathcal{C}_l - \dim S_l = 1$ by a simple dimension count. 
		The result hence follows from the miracle flatness theorem, cf. \cite[Exercise III.10.9]{HartsBook}.
	\end{proof}
	\begin{cor}\label{cor_welldef_slo}
		Assume $c_1(F)$ is very ample.
		Then for any torsion-free coherent sheaf $\mathscr{E}$ over $B$, $l \in \nat^{*(m - 1)}$, there is a proper analytic subset $A_l \subset U_l$, so that the Harder-Narasimhan slopes of the pull-back, $\iota^* \mathscr{E}$, are the same for an arbitrary regular curve $\iota: C \hookrightarrow B$, given by a complete intersection in $l c_1(F)$ and corresponding in $S_l$ to a point away from $A_l$.
	\end{cor}
	\begin{proof}
		Remark first that as $\mathscr{E}$ is torsion-free, it is locally free in codimension at least $2$.
		So over general curves $\iota: C \hookrightarrow B$ given by a complete intersection in $l c_1(F)$, the pull-back $\iota^* \mathscr{E}$ is locally free.
		From now, Corollary \ref{cor_welldef_slo} is a direct consequence of Propositions \ref{prop_cnst_generic}, \ref{prop_bertini}, \ref{prop_flat_fam_inc}.
	\end{proof}
	\par 
	In the following, we will consider several vector bundles and compare the Harder-Narasimhan filtrations of their pull-backs to general curves at a given fiber. 
	The proposition below shows that, to compute the Harder-Narasimhan slopes associated with these pull-backs, it suffices to consider curves passing through a fixed point.
	More precisely, we will need the following result.
	\begin{prop}\label{prop_exist_b}
		Let $\mathscr{E}_k$, $k \in \nat$, be a sequence of torsion-free coherent sheaves over $B$ and $\Sigma \subset B$ be a countable union of proper analytic subsets.
		Then there is $b \in B \setminus \Sigma$ such that for any $l \in \nat^{*(m - 1)}$, there is a curve $\iota: C \hookrightarrow B$, given by a complete intersection in $l c_1(F)$ and passing through $b$, such that for any $k \in \nat$, the pull-back of $\mathscr{E}_k$ to $C$ is torsion-free, and the Harder-Narasimhan slopes of $\iota^* \mathscr{E}_k$ coincide with the Harder-Narasimhan slopes of the pull-back of $\mathscr{E}_k$ to a general curve $C_l \subset B$ given by a complete intersection in $l c_1(F)$.
	\end{prop}
	\par 
	In order to establish Proposition \ref{prop_exist_b}, for a given $l \in \nat^{*(m - 1)}$, we define the subset $W_l \subset B$ in such a way that for any $b \in W_l$, there is a curve $C \subset B$, as required in Proposition \ref{prop_exist_b}.
	\begin{lem}\label{prop_w_l_subs}
		The subset $W_l \subset B$ is a countable intersection of open dense subsets in $B$ (in the manifolds topology). 
	\end{lem}
	\begin{proof}
		We denote by $A_{k, l} \subset U_l$ the subset constructed as $A_l \subset U_l$ from Corollary \ref{cor_welldef_slo} for $\mathscr{E} := \mathscr{E}_k$.
		Remark that $W_l$ has the following description 
		\begin{equation}\label{eq_wl_form}
			W_l = p_{0, l}(\mathcal{C}_l \setminus \cup_{k = 0}^{+ \infty} p_l^{-1}(A_{k, l})).
		\end{equation}				
		It is clearly sufficient to prove that $W_l$ is a dense $G_{\delta}$ subset (recall that $G_{\delta}$ subsets are those which can be written as a countable intersection of open sets with respect to the manifolds topology).
		\par 
		If $W_l$ was not dense, then there would be a non-empty open subset $U \subset B$, such that $W_l \cap U = \emptyset$.
		This would mean that the open subset $p_{0, l}^{-1}(U)$ would lie inside of $\cup_{k = 0}^{+ \infty} p_l^{-1}(A_{k, l})$, which is absurd, as $p_l^{-1}(A_{k, l})$ are proper analytic subsets (in particular, they are nowhere dense and closed), and so their union is nowhere dense by the Baire category theorem.
		\par 
		It is now only left to establish that $W_l$ is $G_{\delta}$. 
		Note that $\mathcal{C}_l \setminus \cup_{k = 0}^{+ \infty} p_l^{-1}(A_{k, l})$ is obviously $G_{\delta}$, as it can be written as $\cap_{k = 0}^{+\infty} V_k$, where $V_k := \mathcal{C}_l \setminus \cup_{i = 0}^{k} p_l^{-1}(A_{i, l})$.
		We claim that 
		\begin{equation}\label{eq_incl_p0l}
			p_{0, l}(\mathcal{C}_l \setminus \cup_{k = 0}^{+ \infty} p_l^{-1}(A_{k, l})) 
			=
			\cap_{k = 0}^{+\infty} p_{0, l}(V_k). 
		\end{equation}
		Note that once it is done, it will follow that $W_l$ is $G_{\delta}$ by (\ref{eq_wl_form}), as $p_{0, l}$ being a locally trivial bundle, is a submersion (and hence open).
		\par 
		Note that the inclusion  $p_{0, l}(\mathcal{C}_l \setminus \cup_{k = 0}^{+ \infty} p_l^{-1}(A_{k, l})) 
			\subset
			\cap_{k = 0}^{+\infty} p_{0, l}(V_k)$ is immediate.
		To establish an inverse inclusion, we pick $b \not\in p_{0, l}(\mathcal{C}_l \setminus \cup_{k = 0}^{+ \infty} p_l^{-1}(A_{k, l}))$ and assume for the sake of a contradiction that for any $k \in \nat$, we have $b \in p_{0, l}(V_k)$.
		The condition $b \not\in p_{0, l}(\mathcal{C}_l \setminus \cup_{k = 0}^{+ \infty} p_l^{-1}(A_{k, l}))$ means exactly that $p_{0, l}^{-1}(b)$ lies inside of $\cup_{k = 0}^{+ \infty} p_l^{-1}(A_{k, l})$, and the condition $b \in p_{0, l}(V_k)$ means that $p_{0, l}^{-1}(b)$ does not lie inside of $\cup_{i = 0}^{k} p_l^{-1}(A_{i, l})$ for any $k \in \nat$.
		Note that $p_{0, l}^{-1}(b)$ is irreducible (as an algebraic variety) and since all $p_l^{-1}(A_{k, l})$ are analytic, they either contain $p_{0, l}^{-1}(b)$ or intersect it over a proper subvariety.
		As proper subvarieties are nowhere dense in irreducible varieties, we conclude that $\cup_{k = 0}^{+ \infty} p_l^{-1}(A_{k, l})$ forms a countable union of nowhere dense closed subsets covering $p_{0, l}^{-1}(b)$.
		This contradicts the Baire category theorem, and this finishes the proof of (\ref{eq_incl_p0l}).
	\end{proof}
	\begin{proof}[Proof of Proposition \ref{prop_exist_b}]
		Clearly, any $b \in \cap_{l \in \nat^{*(m - 1)}} W_l \cap (B \setminus \Sigma)$ would satisfy the requirements of Proposition \ref{prop_exist_b}. 
		So it suffices to establish that $\cap_{l \in \nat^{*(m - 1)}} W_l \cap (B \setminus \Sigma) \neq \emptyset$.
		This is a direct consequence of Lemma \ref{prop_w_l_subs}, the fact that $B \setminus \Sigma$ is a countable intersection of open dense subsets and the Baire category theorem.
	\end{proof}
		
	\section{Dependence of the slopes on the degree of a general curve}\label{sect_deg_fam}
	The main goal of this section is to study how the Harder-Narasimhan slopes of pull-backs of a torsion-free coherent sheaf to general curves depend on the degree of the curve.
	For this, we fix a torsion-free coherent sheaf $\mathscr{E}$ of rank $r$ over $B$. 
 	Assume that the multipolarization $[\omega_B]$ is very ample and integral.
	For any $l \in \nat^{*(m - 1)}$, we denote by $\mu_i(\iota^* \mathscr{E}, l)$, $i = 1, \ldots, r$, the Harder-Narasimhan slopes of the pull-backs $\iota^* \mathscr{E}$ to a general curve $\iota: C \hookrightarrow B$, given by a complete intersection in $l [\omega_B]$, cf. Corollary \ref{cor_welldef_slo}.
	We denote by $\mu_{\min}(\iota^* \mathscr{E}, l)$ and $\mu_{\max}(\iota^* \mathscr{E}, l)$ the minimal and maximal slopes respectively.
	We define the sequence of probability measures, $\eta_l$, $l \in \nat^{* (m - 1)}$, on $\real$ as follows
	\begin{equation}\label{eq_seq_prob_nu_l}
		\eta_l := \frac{1}{r} \sum_{i = 1}^{r} \delta \Big[ \frac{\mu_i(\iota^* \mathscr{E}, l)}{l_1 \cdots l_{m - 1}} \Big].
	\end{equation}
	\begin{sloppypar}
	We say that a sequence $a_l \in \real$, $l = (l_1, \ldots, l_{m - 1}) \in \nat^{*(m - 1)}$, is coordinatewise subadditive, i.e. for any $l' \in \nat^{*(m - 2)}$, $l_1^1, l_1^2 \in \nat^*$, we have 
	\begin{equation}
		a_{l^0} \leq a_{l^1} + a_{l^2},
	\end{equation}
	for $l^0 := (l_1^1 + l_1^2, l')$, $l^1 := (l_1^1, l'), l^2 := (l_1^2, l')$, and similarly for other coordinates.
	\begin{thm}\label{thm_coord_subm}
		For any convex function $\phi : \real \to \real$, the sequence of real numbers 
		\begin{equation}
			a_l := l_1 \cdots l_{m - 1} \cdot \int_{\real} \phi(x) d \eta_l(x),
		\end{equation}
		$l = (l_1, \ldots, l_{m - 1}) \in \nat^{*(m - 1)}$, is coordinatewise subadditive.
	\end{thm}
	\end{sloppypar}
	\par 
	Our proof of Theorem \ref{thm_coord_subm} is based on the result of Shatz \cite{Shatz}, which studies Harder-Narasimhan slopes in the family setting and makes a relation between the Harder-Narasimhan slopes of the general fiber of the family and the Harder-Narasimhan slopes of the specialization of the general fiber.
	For this, we need to introduce a certain order on filtrations.
	\par 
	Let $\mathscr{E}$ be a torsion-free coherent sheaf of rank $r$ over a complex Kähler manifold $B$ with a fixed multipolarization $[\omega_B]$.
	We fix an arbitrary filtration of $\mathscr{E}$  by coherent subsheaves 
	$
		\mathscr{E} = \mathcal{F}_h \supset \mathcal{F}_{h - 1} \supset \cdots \supset \mathcal{F}_1 \supset \mathcal{F}_0 = \{0\},
	$
	and define the \textit{graph of the filtration} as the graph of a piecewise linear function, defined on $[0, r]$, as a linear interpolation of nodes $(\rk{\mathcal{F}_i}, \deg(\mathcal{F}_i))$, $i = 0, \ldots, h$, where the degree is calculated with respect to $[\omega_B]$.
	We say that a graph \textit{dominates} another graph if it lies above it.
	\par 
	The following alternative description of this order will be particularly useful in what follows.
	We define the slopes of the filtration as the non-increasing sequence of numbers $\mu_1^{\mathcal{F}}, \ldots, \mu_r^{\mathcal{F}}$, such that the successive slopes $\mu(\mathcal{F}_i / \mathcal{F}_{i - 1})$ appears among $\mu_1^{\mathcal{F}}, \ldots, \mu_r^{\mathcal{F}}$ exactly $\rk{\mathcal{F}_{i} / \mathcal{F}_{i - 1}}$ times for any $i = 1, \ldots, r$.
	We fix two filtrations of $\mathscr{E}$ by subsheaves the successive slopes of which are non-increasing, and denote their slopes by $\mu_1^{\mathcal{F}, j}, \ldots, \mu_r^{\mathcal{F}, j}$, $j = 1, 2$. 
	Then by \cite[(12.1)]{AtiyahBott}, the first filtration dominates the second one if and only if
	\begin{equation}\label{eq_char_ext1}
		\mu^{\mathcal{F}, 1}_1 + \cdots + \mu^{\mathcal{F}, 1}_k \geq \mu^{\mathcal{F}, 2}_1 + \cdots + \mu^{\mathcal{F}, 2}_k,
	\end{equation}
	for any $k = 1, \ldots, r - 1$. 
	Remark that we always have $\mu^{\mathcal{F}, 1}_1 + \cdots + \mu^{\mathcal{F}, 1}_r = \rk{\mathscr{E}} \cdot \mu(\mathscr{E}) =  \mu^{\mathcal{F}, 2}_1 + \cdots + \mu^{\mathcal{F}, 2}_r$.
	Alternatively, we define the probability measures $\mu^{\mathcal{F}, j}$, $j = 1, 2$, as in (\ref{eq_eta_defn}):
	\begin{equation}
		\mu^{\mathcal{F}, j} := \frac{1}{r} \sum_{i = 1}^{r} \delta \big[ \mu^{\mathcal{F}, j}_i \big].
	\end{equation}
	Then by \cite{HornStoch}, cf. \cite[(12.2)]{AtiyahBott}, the above partial order is equivalent to the following one: for any convex function $\phi : \real \to \real$, we have
	\begin{equation}\label{eq_char_ext2}
		\int \phi d \mu^{\mathcal{F}, 1} \geq \int \phi d \mu^{\mathcal{F}, 2}.
	\end{equation}
	More generally, for arbitrary probability measures $\mu_1, \mu_2$ on $\real$, we say that $\mu_1$ dominates $\mu_2$ if the analogue of (\ref{eq_char_ext2}) holds.
	\par 
	The relevance of this partial order to Harder-Narasimhan filtrations comes from the following extremal characterization of them.
	\begin{thm}[{Shatz \cite[Theorem 2 and Remark on p.174]{Shatz}}]\label{thm_shatz}
		The Harder-Narasimhan filtration of a torsion-free coherent sheaf dominates all the other filtrations by coherent subsheaves.
	\end{thm}
	\begin{rem}
		The article \cite{Shatz} is written for symmetric multipolarization, but it readily extends to the general case, as it only relies on Lemma \ref{lem_slope_zero_morph} and on the existence of the relative Harder-Narasimhan filtration, which hold for non-symmetric multipolarizations, cf. Theorem \ref{thm_rel_hn}.
	\end{rem}
	\par 
	As we shall see below, the above partial order is also useful in the study of how Harder-Narasimhan filtrations behave under specialization.
	For this, we assume from now on that our family of manifolds, $p : B \to S$, is defined over a unit disc, i.e. $S = \mathbb{D} := \{z \in \comp : |z| < 1 \}$. 
	We then decompose the central fiber, $i_0 : B_0 \to B$, as $\sum n_j \cdot B_0^j$, where $i_0^j : B_0^j \to B$, $j = 1, \ldots, t$, are irreducible components, which we assume to be smooth for simplicity, and $n_j$ are the multiplicities.
	We define the probability measure $\mu_0$ on $\real$ as follows
	\begin{equation}\label{eq_meas_cent_fiber}
		\mu_0
		=
		\frac{1}{r} \sum_{i = 1}^{r} \delta \Big[ \sum_{j = 1}^{t} n_j \cdot \mu^{0, j}_i \Big],
 	\end{equation}
 	where $\mu^{0, j}_i$, $j = 1, \ldots, t$, $i = 1, \ldots, r$, are the slopes of the Harder-Narasimhan filtrations of $(i_0^j)^* \mathscr{E}$ (which we assume to be torsion-free).
 	We denote by $\mu_*$ the probability measure composed of the slopes of the Harder-Narasimhan filtrations of the pull-back of $\mathscr{E}$ to general fibers (this doesn't depend on the choice of general fiber by Proposition \ref{prop_cnst_generic}).
 	In other words, by Theorem \ref{thm_rel_hn},
 		\begin{equation}\label{eq_meas_cent_fiber0}
 			\mu_* 
 			= 
 			\frac{1}{r} \sum_{i = 1}^{r} \delta \big[ \mu^{*}_i \big],
 		\end{equation}
 		where $\mu^{*}_i$ are the slopes of the pull-back of the relative Harder-Narasimhan filtration $HN(\mathscr{E})_i$, $i = 0, \ldots, h$, from Theorem \ref{thm_rel_hn} to fibers.
 	\begin{thm}\label{thm_special}
 		The probability measure $\mu_0$ dominates $\mu_*$.
 	\end{thm}
 	\begin{rem}
 		When $B = B_0 \times \mathbb{D}$, and $p : B \to \mathbb{D}$ is the natural projection, the result is due to Shatz \cite[Theorem 3]{Shatz}. When $p$ is of relative dimension $1$, see \cite[Proposition 4.3]{MehtaRamMathAnn} for another related result. 
 		Our proof of the general case is very much inspired by \cite{Shatz}. 
 	\end{rem}
 	In the proof of Theorem \ref{thm_special}, we will use the following simple result.
 	\begin{lem}\label{lem_sum_meas}
 		Let $\mu_i^j = \frac{1}{r} \sum_{l = 1}^{r} \delta[a_{i, l}^j]$, $i, j = 1, 2$, be probability measures on $\real$ for some $a_{i, l}^j \in \real$, ordered in such a way that $a_{i, l}^j$ are non-increasing in $l$, and such that $\sum a_{i, l}^1 = \sum a_{i, l}^2$.
 		Assume that the measure $\mu_i^1$ dominates $\mu_i^2$ for $i = 1, 2$. 
 		Then for the measures $\mu^j := \frac{1}{r} \sum_{l = 1}^{r} \delta[a_{1, l}^j + a_{2, l}^j]$, $\mu^1$ dominates $\mu^2$.
 	\end{lem}
 	\begin{proof}
 		It follows directly from the characterization (\ref{eq_char_ext1}).
 	\end{proof}
 	\begin{proof}[Proof of Theorem \ref{thm_special}]
 		By the flatness of the family and the factors $HN(\mathscr{E})_i / HN(\mathscr{E})_{i - 1}$, cf. \cite[Example 20.3.3]{FultonIntTheory}, we conclude that 
 		\begin{equation}\label{eq_meas_cent_fiber1}
 			\mu^{*}_i
 			=
 			\sum_{j = 1}^{t} n_j \cdot \mu \Big( (i_0^j)^*(HN(\mathscr{E})_i / HN(\mathscr{E})_{i - 1} ) \Big).
 		\end{equation}
 		Now, using the notations introduced in (\ref{eq_meas_cent_fiber}), for any $j = 1, \ldots, t$, we define the measures
 		\begin{equation}\label{eq_meas_cent_fiber2}
 			\begin{aligned}
 			&
			\mu_{0, j}
			:=
			\frac{1}{r} \sum_{i = 1}^{r} \delta \big[ \mu^{0, j}_i \big],
			\\
			&
			\mu_{*, j}
			:=
			\frac{1}{r} \sum_{i = 1}^{r} \delta \Big[ \mu \Big( (i_0^j)^*(HN(\mathscr{E})_i / HN(\mathscr{E})_{i - 1} ) \Big) \Big].
			\end{aligned}
 		\end{equation}
 		By Theorem \ref{thm_shatz}, we conclude that the measures $\mu_{0, j}$ dominate $\mu_{*, j}$.
 		We deduce Theorem \ref{thm_special} from this, Lemma \ref{lem_sum_meas} and (\ref{eq_meas_cent_fiber}), (\ref{eq_meas_cent_fiber0}), (\ref{eq_meas_cent_fiber1}), (\ref{eq_meas_cent_fiber2}).
 	\end{proof}
 	\par 
 	We will now apply this to establish Theorem \ref{thm_coord_subm}.
	The following result will be crucial for this.
	Recall that the manifold $S_{l}$ was defined in (\ref{eq_defn_sl_var}).
	\begin{prop}\label{prop_exists_fam}
		For any $l' = (l_2, \ldots, l_{m - 1}) \in \nat^{*(m - 2)}$, $l_1^1, l_1^2 \in \nat^*$, we let $l^1 := (l_1^1, l')$, $l^2 := (l_1^2, l')$, and fix some curves $C_1, C_2 \subset B$, given by a complete intersection in $l^1 [\omega_B]$, $l^2 [\omega_B]$, as a zero set of holomorphic sections $s_1^1 \in H^0(B, F_1^{l_1^1})$, $s_i \in H^0(B, F_i^{l_i})$, $i = 2, \ldots, m - 1$, and $s_1^2 \in H^0(B, F_1^{l_1^2})$, $s_i$, $i = 2, \ldots, m - 1$, respectively.
		Then there is a smooth manifold $\mathcal{C}$ and a flat proper map $p : \mathcal{C} \to \mathbb{D}$, such that the central fiber of $p$ is given by $C_1 + C_2$, and the general fiber is a complete intersection in $l^0 [\omega_B]$, where $l^0 := (l_1^1 + l_1^2, l')$.
		Moreover, for any countable union of analytic subsets in $S_{l^0}$, we can choose $p$ as a base change along a holomorphic disc in $S_{l^0}$, which doesn't lie within any of the analytic subsets.
	\end{prop} 
	\begin{proof}
		Remark first that there is a canonical embedding $\mathbb{P}(H^0(B, F^{l_1^1})) \times \mathbb{P}(H^0(B, F^{l_1^2}))  \hookrightarrow \mathbb{P}(H^0(B, F^{l_1^1 + l_1^2}))$, given by a composition of the Segre map and the multiplication in the section ring of $B$.
		The needed family can then be constructed by the base change of the family $\mathcal{C}_{l^0}$ (see (\ref{eq_family_cl}) for the definition) along a holomorphic disc in $S_{l^0}$, which is not contained in the countable union of analytic subsets from the statement, and which passes through the image of the point representing $(C_1, C_2)$ in $(S_{l^1}, S_{l^2})$, viewed as a point in $S_{l^0}$ through the above embedding.
		\par 
		The flatness of the resulting map is a consequence of Proposition \ref{prop_flat_fam_inc} and the fact that flatness is preserved by base changes, see \cite[Proposition III.9.2]{HartsBook}.
		The identification of the central fiber follows from the nature of the Segre map: if the holomorphic sections $s_1^1, s_2, \ldots, s_{m - 1}$ (resp. $s_1^2, s_2, \ldots, s_{m - 1}$) vanish along $C_1$ (resp. $C_2$) with multiplicity $1$, then the holomorphic sections $s_1^1 \cdot s_1^2, s_2, \ldots, s_{m - 1}$ vanish along $C_1 + C_2$ with multiplicity $1$.
	\end{proof}
	\par 
	We will also need the following technical lemma. 
	\begin{lem}\label{lem_choice}
		Let $B_1, B_2, B$ be compact complex manifolds, and $A_1 \subset B_1 \times B$, $A_2 \subset B_2 \times B$, be proper analytic subsets. 
		Then there are $x_1 \in B_1, x_2 \in B_2, y \in B$, such that $(x_1, y) \notin A_1$ and $(x_2, y) \notin A_2$.
	\end{lem}
	\begin{proof}
		Let $A$ be a proper analytic subspace in $C \times B$, where $C$ is a compact complex manifold. 
		We denote by $D \subset B$, the locus of $x \in B$, such that $A$ doesn't contain $C \times \{x\}$. 
		By using the properness of the projection map $p : A \to B$ and compactness of $A$, we argue similarly to (\ref{eq_dim_max_irr_comp}), that $D$ is open in $B$ (in Zariski topology).
		It is also non-empty since $A$ is a proper subset of $C \times B$.
		By applying this for $A_1$ and $A_2$, we see that the locus of $y \in B$, such that $A_i$ doesn't contain $B_i \times \{y\}$, $i = 1, 2$, is non-empty.
		It suffices then to pick any $y \in B$ from this locus, and choose $x_1, x_2$ arbitrarily away from $A_1, A_2$.
	\end{proof}
	\par 
	\begin{proof}[Proof of Theorem \ref{thm_coord_subm}.]	
		The coordinatewise subadditivity will be established as a consequence of Theorem \ref{thm_special} applied for a specific degeneration.
		\par 
		More precisely, from Corollary \ref{cor_welldef_slo}, the slopes of the Harder-Narasimhan filtrations of $\iota^* \mathscr{E}$ are constant for curves $\iota: C \hookrightarrow B$ given by a complete intersection in $l^i [\omega_B]$, $i = 0, 1, 2$, located away from a proper analytic subset, $A_{l^i}$, in $S_{l^i}$.
		Using Lemma \ref{lem_choice}, we fix curves $C_i \in S_{l^i} \setminus A_{l^i}$, given by complete intersection in $l^i [\omega_B]$, $i = 1, 2$, as required in Proposition \ref{prop_exists_fam}.
		By Proposition \ref{prop_exists_fam}, $C_1 + C_2$ can be put as a central fiber in a flat family of curves, obtained as a base change along a holomorphic disc from $S_{l^0}$, not contained in $A_{l^0}$.
		Since the general curve from this family doesn't lie in $A_{l^0}$, the Harder-Narasimhan slopes of the pull-back of $\mathscr{E}$ to a general curve in this family will then coincide with the Harder-Narasimhan slopes of the pull-back of $\mathscr{E}$ to general curves given by complete intersection in $l^0 [\omega_B]$.
		Directly from this, Theorem \ref{thm_special} and (\ref{eq_char_ext2}),
		\begin{equation}
			\sum_{i = 1}^{r} \phi \big( \mu_i(\iota^* \mathscr{E}, l^1) + \mu_i(\iota^* \mathscr{E}, l^2) \big)
			\geq
			\sum_{i = 1}^{r} \phi \big( \mu_i(\iota^* \mathscr{E}, l^0) \big).
		\end{equation}
		By making a homothety and using the convexity of $\phi$, we finish the proof.
	\end{proof}
	As a byproduct of Theorem \ref{thm_coord_subm} and (\ref{eq_char_ext1}), we deduce the following.
	\begin{prop}\label{prop_unif_bnd_hn_slopes_curves2}
		For any $l = (l_1, \ldots, l_{m - 1}) \in \nat^{*(m - 1)}$, we have
		\begin{equation}
		\begin{aligned}
			&
			\mu_{\min}(\iota^* \mathscr{E}, l)
			\geq 
			l_1 \cdots l_{m - 1} \cdot
			\mu_{\min}(\iota^* \mathscr{E}, (1, \ldots, 1)), 
			\\
			&
			\mu_{\max}(\iota^* \mathscr{E}, l)
			\leq
			l_1 \cdots l_{m - 1} \cdot
			\mu_{\max}(\iota^* \mathscr{E}, (1, \ldots, 1)).
		\end{aligned}
		\end{equation}
	\end{prop}

	\section{General base changes and asymptotic distribution of slopes}\label{sect_MR}
	The main goal of this section is to show that the asymptotic distribution of the Harder-Narasimhan slopes of direct images can be calculated through general base changes, i.e. to establish Theorem \ref{thm_mehta_ramanathan}.
	The following famous result will be at the heart of our approach.
	\begin{thm}\label{thm_mr}
		For any semistable torsion-free coherent sheaf $\mathscr{E}$ over a complex projective manifold $B$ with an integral multipolarization $[\omega_B]$, there is $l_0 \in \nat$, such that for any general curve $\iota: C \hookrightarrow B$, which is a complete intersection in $l [\omega_B]$, $l = (l_1, \ldots, l_{m - 1}) \in \nat^{*(m - 1)}$, $l_i \geq l_0$, $i = 1, \ldots, m - 1$, the pull-back $\iota^* \mathscr{E}$ is semistable.
	\end{thm}
	\begin{rem}
		When the multipolarization is symmetric, the theorem is due to Mehta-Ramanathan \cite[Theorem 6.1]{MehtaRamMathAnn}, cf. also \cite[Theorem 7.2.8]{HuybLehn} and \cite{MehtaRama}.
		The proof of the general multipolarization case follows directly from a result of Langer \cite[Corollary 5.4]{LangerAnnals}, which states that in the multipolarization setting, the pull-back of a semistable sheaf to a general hypersurface associated with the first component of the multipolarization remains semistable if the degree of the hypersurface is big enough. 
		By iterating this statement, one obtains the corresponding result for curves.
	\end{rem}
	In particular, from Theorem \ref{thm_mr} and (\ref{eq_slope_rest}), we conclude that for any torsion-free coherent sheaf $\mathscr{E}$ over $B$, the pull-back of the Harder-Narasimhan filtration of $\mathscr{E}$ over general curves $\iota: C \hookrightarrow B$ as in Theorem \ref{thm_mr} coincides with the Harder-Narasimhan filtration of $\iota^* \mathscr{E}$.
	\par 
	The main difficulty in the application of Theorem \ref{thm_mr} to the proof of Theorem \ref{thm_mehta_ramanathan} is that the value $l_0$ from Theorem \ref{thm_mr} depends heavily on $\mathscr{E}$.
	Despite many works on the effective version of Theorem \ref{thm_mr}, see Flenner \cite{Flenner} or Bogomolov \cite{BogomStable}, cf. \cite[\S 7]{HuybLehn}, it seems that all known effective bounds on $l_0$ depend on the rank and the discriminant of $\mathscr{E}$. 
	However, in the notations from Introduction, when $k \to \infty$, both rank and discriminant of $\mathscr{E}_k$ might tend to infinity.
	\par 
	Remark also that there are some results that relate the Harder-Narasimhan slopes of a sheaf with the Harder-Narasimhan slopes of its pull-back to a curve associated with any fixed value $l_0$, see for example the Grauert-M\"{u}lich theorem, cf. \cite[\S 3]{HuybLehn}, or Langer \cite[Theorem 3.1]{LangerAnnals}.
	To the best of the author's knowledge, none of these implies Theorem \ref{thm_mehta_ramanathan} directly. 
	Our approach for Theorem \ref{thm_mehta_ramanathan} is then, by necessity, a different one, and it relies on the subadditivity properties enjoyed by the measures $\eta_l^{HN|C}$ and $\eta^{HN}$ established in Section \ref{sect_deg_fam}.
	\par 
	We first establish the identity (\ref{eq_mehta_ramanathan_minmax0}) about the minimal slopes.
	For this, we conserve the notations from Introduction and (\ref{eq_seq_prob_nu_l}).
	Directly from Theorem \ref{thm_mr}, the remark below it and Lemma \ref{lem_major_hn_restr}, in the notations analogous to (\ref{eq_eta_cmax_edefe}), we deduce the following result.	
	\begin{lem}\label{lem_min_slope_conv}
		The following convergence takes place
		\begin{equation}\label{eq_min_slope_conv1}
			\lim_{l \to \infty} \eta_{\min, k, l}^{HN|C}
			=
			\frac{\mu_{\min}(\mathscr{E}_k)}{k},
		\end{equation}
		where $l = (l_1, \ldots, l_{m - 1}) \in \nat^{* (m - 1)}$.
		Moreover, the above limit becomes stationary for $l$ big enough (i.e. for big enough $l_0 \in \nat$, and $l_i \geq l_0$, $i = 1, \ldots, m - 1$), and we always have
		\begin{equation}\label{eq_min_slope_conv2}
			\eta_{\min, k, l}^{HN|C}
			\leq
			\frac{\mu_{\min}(\mathscr{E}_k)}{k}.
		\end{equation}
	\end{lem}
	\begin{rem}\label{rem_max_sl_inv}
		The same proof shows that the analogous statements hold (with an inverse sign in (\ref{eq_min_slope_conv2})) for the maximal slope.
	\end{rem}
	We will also use the following statement, the proof of which is left to the reader.
	\begin{lem}\label{lem_lim_ord_0}
		Let $c_{k, l} \in \real$, $k \in \nat$, $l \in \nat^{* (m - 1)}$, be such that $c_{k, l}$ increases in $k \in \nat$, and the limit $\lim_{l \to \infty} c_{k, l}$ exists and satisfies $c_{k, l} \leq \lim_{l \to \infty} c_{k, l}$.
		Then the double limits $\lim_{l \to \infty} \lim_{k \to \infty} c_{k, l}$, $\lim_{k \to \infty} \lim_{l \to \infty} c_{k, l}$ exist and they are equal.
	\end{lem}
	\begin{sloppypar}
	\begin{proof}[Proof of (\ref{eq_mehta_ramanathan_minmax0})]
		By (\ref{eq_min_slope_conv1}) and the definitions of $\eta_{\min, l}^{HN|C}$ and $\eta_{\min}^{HN}$, we see that (\ref{eq_mehta_ramanathan_minmax0}) is equivalent to the following statement
		\begin{equation}\label{eqqlem_lim_ord_0}
			\lim_{l \to \infty} \lim_{k \to \infty} \eta_{\min, k, l}^{HN|C}
			=
			\lim_{k \to \infty} \lim_{l \to \infty} \eta_{\min, k, l}^{HN|C}.
		\end{equation}
		In fact, by the eventual superadditivity of the minimal slopes, see Proposition \ref{prop_superadd} and \ref{prop_sm}, it suffices to establish the above statement when the limit over $k$ runs over a subsequence of $\nat$, as the limits exist by Fekete's lemma.
		From now on, we assume that $k = 2^r$, $r \in \nat$.
		Then by Proposition \ref{prop_superadd}, the sequence $\eta_{\min, k, l}^{HN|C}$ increases in $k$. 
		By this and Lemma \ref{lem_min_slope_conv}, the assumptions of Lemma \ref{lem_lim_ord_0} are satisfied for $c_{k, l} := \eta_{\min, k, l}^{HN|C}$.
		Applying Lemma \ref{lem_lim_ord_0} then yields (\ref{eqqlem_lim_ord_0}).
	\end{proof}
	\end{sloppypar}
	We will now concentrate on the proof of Theorem \ref{thm_mehta_ramanathan}.
	We pick $b \in B$ as in Proposition \ref{prop_exist_b}, applied for direct images $\mathscr{E}_k$ and the set $\Sigma := \Sigma \cup \Sigma'$ as in the beginning of Section \ref{sect_sm_bnd} and Proposition \ref{prop_sm}.
	Let $i_b : X_b \to X$ be the fiber of $\pi$ at $b$ and $L_b := i_b^* L$.
	We fix $l \in \nat^{* (m - 1)}$ and consider a general curve $\iota: C \hookrightarrow B$ given by a complete intersection in $l [\omega_B]$ passing through $b$ as in Proposition \ref{prop_exist_b}.
	We denote by $\mathcal{F}^C_b(l)$ the filtration on $R(X_b, L_b)$ induced by (\ref{eq_map_restriction_sheaf}) and the Harder-Narasimhan filtration on $\iota^* \mathscr{E}_k$.
	For any $d \in \nat^*$, we consider the induced filtration $\mathcal{F}^{C, d}_b(l)$ on the subalgebra $[{\rm{Sym}} H^0(X_{b}, L_b^d)]$. 
	We denote by $\eta_{k, l}^{d, 0}$ the probability measure on $\real$, defined as
	\begin{equation}\label{eq_ind_measu_defn_HN}
		\eta_{k, l}^{d, 0}
		=
		\frac{1}{\dim [{\rm{Sym}}^k H^0(X_{b}, L_b^d)]}
		\sum_{j = 1}^{\dim [{\rm{Sym}}^k H^0(X_{b}, L_b^d)]}
		\delta \Big[ 
			\frac{\eta_{\mathcal{F}^{C, d}_b(l)}(j, k)}{k d \cdot l_1 \cdots l_{m - 1}}
		\Big].
	\end{equation}
	Let $\eta_{l}^{d, 0}$ be the weak limit of $\eta_{k, l}^{d, 0}$, as $k \to \infty$ (it exists by the results of Section \ref{sect_gen_th_sf}).
	We denote by $\eta_{k}^{d}$ probability measure on $\real$, defined as the jumping measure of the filtration on $[{\rm{Sym}}^k H^0(X_{b}, L_b^d)]$ induced by the restriction of the Harder-Narasimhan filtration on $\mathscr{E}_d$, i.e. by the procedure outlined before Theorem \ref{thm_filt}.
	We denote by $\eta^{d}$ the associated limiting jumping measure on $\real$.
	\par 
	Our proof of Theorem \ref{thm_mehta_ramanathan} will be based on a careful comparison of these measures.
	As a first step in this direction, remark that as the Harder-Narasimhan filtration is submultiplicative by Proposition \ref{prop_sm2}, and the filtration leading to the definition of $\eta_{k}^{d}$ is the maximal submultiplicative filtration on $[{\rm{Sym}}^k H^0(X_{b}, L_b^d)]$  among those which coincide with the restriction of the Harder-Narasimhan filtration on $\mathscr{E}_d$ on $H^0(X_{b}, L_b^d)$, by Theorem \ref{thm_mr}, for any increasing function $\phi : \real \to \real$ and any $d \in \nat^*$, there is $l_0 \in \nat$, so that for any $k \in \nat$, $l \geq l_0$, we have
	\begin{equation}\label{eq_monoton_contr}
		\int \phi d \eta_{k}^{d}
		\geq
		\int \phi  d\eta_{k, l}^{d, 0}.
	\end{equation}
	We will also use the following easy statement.
	\begin{lem}\label{lem_int_phi_m_comp}
		Assume that a continuous function $\phi : \real \to \real$ has a bounded from below support. 
		Then for any $\epsilon > 0$, there is $d_0 \in \nat$, so that for any $d \geq d_0$, there is $k_0 \in \nat$, so that for any $k \geq k_0$ and $l \in \nat^{* (m - 1)}$, we have
		\begin{equation}
			\Big|
				\int \phi d \eta_{kd, l}^{HN|C}
				-
				\int \phi d \eta_{k, l}^{d, 0}
			\Big| \leq \epsilon.
		\end{equation}
	\end{lem}
	Lemma \ref{lem_int_phi_m_comp} follows immediately from Theorem \ref{thm_sym_subalgebra} and the following two results.
	\begin{lem}[{ Chen \cite[Proposition 1.2.5]{ChenHNolyg} }]\label{lem_shr_exact_seq}
		Let $0 \rightarrow V_0 \rightarrow V_1 \rightarrow V_2 \rightarrow 0$ be a short exact sequence of vector spaces.
		We fix a filtration $\mathcal{F}_1$ on $V_1$, and induce the filtrations $\mathcal{F}_0$, $\mathcal{F}_2$ on $V_0$ and $V_2$, which we denote by an abuse of notation by the same letter.
		Then the following identity holds
		\begin{equation}
			\dim V_1 \cdot \eta_{\mathcal{F}_1} = \dim V_0 \cdot \eta_{\mathcal{F}_0} + \dim V_2 \cdot \eta_{\mathcal{F}_2},
		\end{equation}
		where the jumping measures, $\eta_{\mathcal{F}_i}$, $i = 0, 1, 2$, are defined in (\ref{eq_jump_meas_gen}).
	\end{lem}
	\begin{prop}\label{prop_unif_bnd_hn_slopes_curves}
		For any $d \in \nat^*$, there is $C > 0$, such that for any $k \in \nat^*$, $l \in \nat^{* (m - 1)}$,
		\begin{equation}
			\esssup \eta_{k, l}^{d, 0}
			\leq
			C.
		\end{equation}
	\end{prop}
	\begin{proof}
		It follows from Propositions \ref{thm_bound2} and \ref{prop_unif_bnd_hn_slopes_curves2}.
	\end{proof}
	\begin{proof}[Proof of Theorem \ref{thm_mehta_ramanathan}]
		We first establish the vague convergence of measures. 
		As we shall explain later, the statement about the convergence of maximal slopes will then follow from it.
		\par 
		\begin{sloppypar}
		To verify vague convergence, it suffices to verify that for any function $\phi : \real \to \real$ of the form $\phi(x) = \max\{0, x - \lambda \}$ for a certain $\lambda \in \real$, as $l \to \infty$, we have
		\begin{equation}\label{eq_thm_mehta_ramanathan_1}
			\lim_{l \to \infty} \lim_{k \to \infty} \int \phi d \eta_{k, l}^{HN|C}
			=
			\int \phi d \eta^{HN}.
		\end{equation}
		Indeed, such functions determine the vague convergence as any smooth function of compact support $g : \real \to \real$ writes as $g(x) = \int_{-\infty}^x \max\{0, x - t \} \cdot g''(t) dt$.
		\par 
		Let us first establish the following bound
		\begin{equation}\label{eq_conv_two_sid_first}
			\liminf_{l \to \infty} \lim_{k \to \infty} \int \phi d \eta_{k, l}^{HN|C}
			\geq
			\int \phi d \eta^{HN}.
		\end{equation}
		From Theorem \ref{thm_mr} and the remark after it, we have
		\begin{equation}\label{eq_thm_mehta_ramanathan_3}
			\int \phi d \eta^{HN}
			=
			\lim_{k \to \infty} \lim_{l \to \infty} \int \phi d \eta_{k, l}^{HN|C}.
		\end{equation}
		Note that by Theorem \ref{thm_coord_subm},  the sequence $(k, l) \mapsto l_1 \cdots l_{m - 1} \int \phi d \eta_{k, l}^{HN|C}$, $k \in \nat^*$, $l = (l_1, \ldots, l_{m - 1}) \in \nat^{*(m - 1)}$, is coordinatewise subadditive in $l$. 
		By Fekete's lemma, $\int \phi d \eta_{k, l}^{HN|C} \geq \lim_{l' \to \infty} \int \phi d \eta_{k, l'}^{HN|C}$.
		By taking limits $k \to \infty$ and then $l \to \infty$, using (\ref{eq_thm_mehta_ramanathan_3}), we obtain (\ref{eq_conv_two_sid_first}).
		\end{sloppypar}
		\par 
		Let us now establish the inverse inequality 
		\begin{equation}\label{eq_sss11}
			\int \phi d \eta^{HN}
			\geq
			\limsup_{l \to \infty} \lim_{k \to \infty} \int \phi d \eta_{k, l}^{HN|C}.
		\end{equation}
		Note that the limit in $k$-variable in (\ref{eq_thm_mehta_ramanathan_3}) exists by Theorem \ref{thm_filt}, and so it suffices to verify (\ref{eq_sss11}) where $k$ runs over a subsequence.
		By this and Lemma \ref{lem_int_phi_m_comp}, it suffices to show that for any $\epsilon > 0$, there is $d_0 \in \nat^*$, so that for any $d \geq d_0$, there is $l_0 \in \nat$, so that for any $l \geq l_0$, we have
		\begin{equation}\label{eq_thm_mehta_ramanathan_5}
			\int \phi d \eta^{HN}
			\geq
			\int \phi d \eta_{l}^{d, 0} - \epsilon,
		\end{equation}
		on the proof of which we concentrate from now on.
		\par 
		Note that by Theorem \ref{thm_filt}, there is $d_0 \in \nat^*$, so that for any $d \geq d_0$, we have
		\begin{equation}\label{eq_vol_bnd_below}
			\Big|
			\int \phi d \eta^{HN}
			-
			\int \phi d \eta^d
			\Big|
			\leq
			\epsilon.
		\end{equation}
		We now fix $d \geq d_0$.
		Then by (\ref{eq_monoton_contr}), there is $l_0 \in \nat$, so that for any $l \geq l_0$, we have
		\begin{equation}\label{eq_vol_bnd_below2}
			\int \phi d \eta^d
			\geq
			\int \phi d \eta_l^{d, 0}.
		\end{equation}
		Combining (\ref{eq_vol_bnd_below}) and (\ref{eq_vol_bnd_below2}), yields (\ref{eq_thm_mehta_ramanathan_5}).
		\par 
		\begin{sloppypar}
		We will now establish the second convergence of (\ref{eq_mehta_ramanathan_minmax}).
		Note that if $\kappa(L_f) = 0$, the result follows immediately from the weak convergence of the measures.
		We will hence assume from now on that $\kappa(L_f) > 0$.
		First of all, from Lemma \ref{lem_major_hn_restr}, we see
		\begin{equation}\label{eq_one_sided_ineq}
			\liminf_{l \to \infty} \eta_{\max, l}^{HN|C} \geq \eta_{\max}^{HN}.
		\end{equation}
		Now, by Theorem \ref{thm_filt}, the cumulative function $F_{HN}(\lambda) := \int_{\lambda}^{\infty} d\eta^{HN}$ (resp. $F_{HN, l}(\lambda) := \int_{\lambda}^{\infty} d \eta_l^{HN|C}$) is such that $F_{HN}^{1 / \kappa(L_f)}$ (resp. $F_{HN, l}^{1 / \kappa(L_f)}$) is concave on $]- \infty, \esssup \eta^{HN}[$. In particular, the function $F_{HN}$ is continuous on $\real \setminus \{ \esssup \eta^{HN} \}$. 
		As vague convergence of compactly supported measures implies pointwise convergence of the cumulative functions over the continuity set.
		In particular, $F_{HN, l}(\lambda)$ converge pointwise, as $l \to \infty$, to $F_{HN}(\lambda)$ away from $\esssup \eta^{HN}$. 
		\end{sloppypar}
		\par
		Assume now, for the sake of contradiction, that we have a strict inequality
		\begin{equation}\label{eq_one_sided_ineq2}
			\limsup_{l \to \infty} \eta_{\max, l}^{HN|C} > \eta_{\max}^{HN}.
		\end{equation}
		Upon restricting to a subsequence in $l$, we may assume that we have $\liminf$ instead of $\limsup$ in (\ref{eq_one_sided_ineq2}).
		Let $\epsilon > 0$ be such that the inequality (\ref{eq_one_sided_ineq2}) is still satisfied if we replace $\eta_{\max}^{HN}$ by $\eta_{\max}^{HN} + 4 \epsilon$.
		We assume that $l$ is big enough so that $\eta_{\max}^{HN} + 2 \epsilon < \eta_{\max, l}^{HN|C}$.
		Remark that by Theorem \ref{thm_filt}, $\eta_{\max}^{HN}$ coincides with the essential supremum of $\eta^{HN}$, hence $F_{HN}(\eta_{\max}^{HN} - \epsilon)$ is strictly positive.
		\par 
		Note that by the Portemanteau theorem, $F_{HN, l}(\eta_{\max}^{HN} - \epsilon)$ converges, as $l \to \infty$, to $F_{HN}(\eta_{\max}^{HN} - \epsilon)$.
		Since $F_{HN}(\eta_{\max}^{HN} - \epsilon)$ is strictly positive and $F_{HN, l}(\eta_{\max}^{HN} + 2 \epsilon)$ is non-negative, by the concavity of $F_{HN, l}^{1 / \kappa(L_f)}$ on the interval $]-\infty, \eta_{\max}^{HN} + 2\epsilon] \subset ]-\infty, \eta_{\max, l}^{HN|C}[$, there is $\delta > 0$, such that
		\begin{equation}
			F_{HN, l}(\lambda) > \delta, \qquad \text{for any } \lambda \in ]-\infty, \eta_{\max}^{HN} + \epsilon],
		\end{equation}
		and any $l$ in the fixed subsequence.
		This clearly contradicts the fact that $F_{HN, l}(\eta_{\max}^{HN} + \epsilon)$ converges to $F_{HN}(\eta_{\max}^{HN} + \epsilon) = 0$, as $l \to \infty$.
		Hence, the initial assumption (\ref{eq_one_sided_ineq2}) was false.
	\end{proof}

\section{Multiplicities of line bundles and the Harder-Narasimhan slopes}\label{sect_vol_est}
	The main goal of this section is to show that in a fibered setting, multiplicities of line bundles defined over a family of varieties can be estimated in terms of the induced Harder-Narasimhan filtration and the multiplicities of the fiber. In other words, we establish Theorem \ref{thm_vol_bnd}.
	\par 
	Let us recall the following useful statement. 
	We fix a compact connected Riemann surface $C$ of genus $g$.
	For any holomorphic vector bundle $E$ on $C$, we denote by ${\rm{deg}}_+(E)$ its positive degree, defined as the sum of the non-negative Harder-Narasimhan slopes of $E$.
	\begin{prop}[{cf. \cite[Lemma 2.1]{ChenVolume}}]\label{prop_gl_sect_rs}
		We have
		\begin{equation}\label{eq_gl_sect_rs}
			\Big| \dim H^0(C, E) - {\rm{deg}}_+(E) \Big| \leq (g + 1) \cdot \rk{E}.
		\end{equation}
		Moreover, if the Harder-Narasimhan slopes of $E$ are negative, then $H^0(C, E) = \{0\}$.
	\end{prop}
	We can now proceed to the proof of the main result of the current section.
	\begin{proof}[Proof of Theorem \ref{thm_vol_bnd}]
		Let us first establish the bound (\ref{eq_kappa_ineq}).
		We denote by $Y_k$ the image of $X$ under the Kodaira-Iitaka maps $X \dashrightarrow \mathbb{P}(H^0(X, L^k)^*)$.
		By the projection formula, the following isomorphism takes place
		\begin{equation}\label{eq_projjj}
			H^0(X, L^k)
			\simeq
			H^0(B, \mathscr{E}_k).
		\end{equation}
		In particular, we have a natural map $B \times H^0(X, L^k) \to \mathscr{E}_{k}$, inducing the rational map $\mathbb{P}(\mathscr{E}_{k}^*) \dashrightarrow \mathbb{P}(H^0(X, L^k)^*)$.
		We moreover have the fiberwise Kodaira-Iitaka map $X \dashrightarrow \mathbb{P}(\mathscr{E}_{k}^*)$ which commutes with the projection to $B$.
		We denote by $Y_k'$ the image of this map.
		Since the maps $X \dashrightarrow \mathbb{P}(H^0(X, L^k)^*)$ factorize through $X \dashrightarrow \mathbb{P}(\mathscr{E}_{k}^*)$, we have the induced dominant map $Y_k' \dashrightarrow Y_k$, the existence of which implies that 
		\begin{equation}\label{eq_kod_yoga_1}
			\dim Y_k \leq \dim Y_k'.
		\end{equation}
		On another hand, the fiber at a general $b \in B$ of $Y_k'$ coincides with the image $Y'_{k, b}$ of the fiber $X_b$ of $\pi$ under the respective fiberwise Kodaira-Iitaka map.
		Moreover, by \cite[Exercise II.3.22]{HartsBook}, we deduce that for general $b \in B$, we have
		\begin{equation}\label{eq_kod_yoga_2}
			\dim Y'_k = \dim Y'_{k, b} + \dim B.
		\end{equation}
		The bound (\ref{eq_kappa_ineq}) then follows from (\ref{eq_dim_ki_yk}), (\ref{eq_kod_yoga_1}) and (\ref{eq_kod_yoga_2}).
		\par 
		Let us first establish the estimate (\ref{eq_vol_bnd}).
		For this, for any $l = (l_1, \ldots, l_{m - 1}) \in \nat^{m - 1}$ big enough, we establish that in the notations of Theorems \ref{thm_mehta_ramanathan} and \ref{thm_vol_bnd}, the following bound holds
		\begin{equation}\label{eq_vol_bnd_new}
			\lim_{k \to \infty} \frac{\dim H^0(X, L^{k})}{k^{\kappa(L_f) + m}}
			\leq
			\frac{{\rm{vol}}_{\kappa}(L_f)}{\kappa(L_f)! m! \cdot \prod_{i = 1}^{m - 1} a_i}
			\cdot
			\int_{x = 0}^{+\infty} x^m \cdot d \eta_l^{HN|C}(x).
		\end{equation}
		The estimate (\ref{eq_vol_bnd}) will then follow immediately from Theorem \ref{thm_mehta_ramanathan} and (\ref{eq_vol_bnd_new}).
		\par 
		Let us now proceed to the proof of (\ref{eq_vol_bnd_new}).
		By (\ref{eq_projjj}), it is sufficient to estimate $\dim H^0(B, \mathscr{E}_k)$, from which we concentrate from now on.
		\par 
		We fix $l = (l_1, \ldots, l_{m - 1}) \in \nat^{m - 1}$ big enough, so that for any $i = 1, \ldots, m - 1$, there is a smooth divisor $D_i$ in the class $l_i  [\omega_{B, i}]$, which we assume to be general.
		For $j = 1, \ldots, m - 1$, we then denote $B_j := \cap_{i = 1}^{j} D_i$, $\iota_j : B_j \to B$, and we let $B_0 := B$.
		We also denote by $D'_j$ the divisor $D_j \cap B_{j - 1}$ on $B_{j - 1}$.
		Then for any $j = 1, \ldots, m - 1$, and any coherent sheaf $\mathscr{E}$ on $B$, we have
		\begin{equation}\label{eq_bj_iter_e}
			\dim H^0(B_{j - 1}, iota_{j - 1}^* \mathscr{E})
			\leq
			\sum_{i_j = 0}^{+\infty}
			\dim H^0(B_j, \iota_j^* ( \mathscr{E} \otimes \mathscr{O}(-i_j D'_j))).
		\end{equation}
		Indeed, the above bound is obtained by summing the contributions arising from the long exact cohomology sequence associated with the short exact sequence of sheaves
		\begin{equation}
			0 \to  \mathscr{E} \otimes \mathscr{O}(-(i_j + 1)  D'_j) \to \mathscr{E} \otimes \mathscr{O}(-i_j  D'_j) \to \iota_{j}^* (\mathscr{E} \otimes \mathscr{O}(-i_j  D'_j)) \to 0,
		\end{equation}
		and then applying Serre duality together with the vanishing theorem to assure that the sum above is in fact finite.
		\par 
		We now denote $C := B_{m - 1}$, $\iota : C \hookrightarrow B$, which is a compact (smooth) Riemann surface by Bertini's theorem.
		For any $j = 1, \ldots, m - 1$, we denote by $L_j$ the pull-back of the line bundle $\mathscr{O}(D_j)$ to $C$.
		By iterating the inequalities (\ref{eq_bj_iter_e}), applied for $\mathscr{E}_k$ and their twists by the divisorial line bundles, for any $k \in \nat$, we obtain the following inequality
		\begin{equation}\label{eq_bnd_env_sum_div}
			\dim H^0(B, \mathscr{E}_k)
			\leq
			\sum_{i_1, \cdots, i_{m - 1} = 0}^{+\infty}
			\dim H^0(C, \iota^* \mathscr{E}_k \otimes \otimes_{j = 1}^{m - 1} L_j^{-i_j}),
		\end{equation}
		where the sum above is again in fact finite.
		\par 
		We will now apply Proposition \ref{prop_gl_sect_rs} to estimate every cohomology group above.
		For any $a \in \real$, we use the notation $(x - a)_+ := \max\{ x - a, 0 \}$.
		To simplify our further formulas, we let $\overline{\deg}(L_j) := \deg(L_j)/(l_1 \cdots l_{m - 1})$.
		In the notations introduced before Theorem \ref{thm_mehta_ramanathan}, by Proposition \ref{prop_gl_sect_rs}, there is $c > 0$, so that
		\begin{multline}\label{eq_bnd_env_sum_div1}
			\dim H^0(C, \iota^* \mathscr{E}_k \otimes \otimes_{j = 1}^{m - 1} L_j^{-i_j})
			\\
			\leq
			\rk{\mathscr{E}_k} \cdot \Big( c + k \cdot l_1 \cdots l_{m - 1} \cdot \int_{x \in \real} \big( x - \frac{1}{k} \sum_{j = 1}^{m - 1} i_j \overline{\deg}(L_j) \big)_+ d \eta_{k, l}^{HN|C}(x) \Big).
		\end{multline}
		Note that for any $b \geq 0$, the following bound holds $(x + b)_+ \leq x_+ + b$.
		From this, we obtain
		\begin{multline}\label{eq_bnd_env_sum_div2}
			\int_{x \in \real} \big( x - \frac{1}{k} \sum_{j = 1}^{m - 1} i_j \overline{\deg}(L_j) \big)_+ d \eta_{k, l}^{HN|C}(x)
			\leq
			\frac{c_0 \cdot m}{k}
			\\
			+
			\int_{t_1 = i_1}^{i_1 + 1} \cdots \int_{t_{m - 1} = i_{m - 1}}^{i_{m - 1} + 1} 
			\int_{x \in \real} \big( x - \frac{1}{k} \sum_{j = 1}^{m - 1} t_j  \overline{\deg}(L_j) \big)_+ d \eta_{k, l}^{HN|C}(x) dt_1 \cdots dt_{m - 1},
		\end{multline}
		where $c_0 := \max\{ \overline{\deg}(L_1), \ldots, \overline{\deg}(L_{m - 1}) \}$.
		We can now rewrite
		\begin{multline}\label{eq_bnd_env_sum_div3}
			\sum_{i_1, \cdots, i_{m - 1} = 0}^{+\infty} 
			\int_{t_1 = i_1}^{i_1 + 1} \cdots \int_{t_{m - 1} = i_{m - 1}}^{i_{m - 1} + 1} 
			\int_{x \in \real} \big( x - \frac{1}{k} \sum_{j = 1}^{m - 1} t_j \overline{\deg}(L_j) \big)_+ d \eta_{k, l}^{HN|C}(x) dt_1 \cdots dt_{m - 1}
			\\
			=
			\int_{t_1 = 0}^{+ \infty} \cdots \int_{t_{m - 1} = 0}^{+ \infty} 
			\int_{x \in \real} \big( x - \frac{1}{k} \sum_{j = 1}^{m - 1} t_j \overline{\deg}(L_j) \big)_+ d \eta_{k, l}^{HN|C}(x) dt_1 \cdots dt_{m - 1}.
		\end{multline}
		\par 
		Note also that by Propositions \ref{prop_unif_bnd_hn_slopes_curves2} and the second part of Proposition \ref{prop_gl_sect_rs}, for any $l \in \nat^{m - 1}$, there is $c' \in \nat$, so that in the sum on the right-hand side of (\ref{eq_bnd_env_sum_div}), only the terms corresponding to $i_j =  0, \ldots, c' k$ contribute in a nontrivial way.
		Combining this, (\ref{eq_bnd_env_sum_div}), (\ref{eq_bnd_env_sum_div1}), (\ref{eq_bnd_env_sum_div2}) and (\ref{eq_bnd_env_sum_div3}), we deduce that for any $l \in \nat^{m - 1}$, there is $c'' > 0$, so that for $k \in \nat$ big enough, we have
		\begin{multline}
			\dim H^0(B, \mathscr{E}_k)
			\leq
			\rk{\mathscr{E}_k}
			\cdot l_1 \cdots l_{m - 1}
			\cdot
			\\
			\cdot
			\Big(
			k \cdot \int_{t_1 = 0}^{+ \infty} \cdots \int_{t_{m - 1} = 0}^{+ \infty} 
			\int_{x \in \real} \big( x - \frac{1}{k} \sum_{j = 1}^{m - 1} t_j  \overline{\deg}(L_j) \big)_+ d \eta_{k, l}^{HN|C}(x) dt_1 \cdots dt_{m - 1}
			+
			c'' k^{m - 1}
			\Big).
		\end{multline}
		From the definition of the multiplicity, by (\ref{eq_projjj}), we obtain
		\begin{multline}\label{eq_vol_fin_1}
			\lim_{k \to \infty} \frac{\dim H^0(X, L^{k})}{k^{\kappa(L_f) + m}}
			\leq
			\frac{{\rm{vol}}_{\kappa}(L_f)}{\kappa(L_f)!}	
			\cdot l_1 \cdots l_{m - 1}
			\cdot
			\\
			\cdot
			\limsup_{k \to \infty}
			\frac{1}{k^{m - 1}}
			\int_{t_1 = 0}^{+ \infty} \cdots \int_{t_{m - 1} = 0}^{+ \infty} 
			\int_{x \in \real} \big( x - \frac{1}{k} \sum_{j = 1}^{m - 1} t_j  \overline{\deg}(L_j) \big)_+ d \eta_{k, l}^{HN|C}(x) dt_1 \cdots dt_{m - 1}.
		\end{multline}
		\par 
		By making a homothety for every variable $t_j$, $j = 1, \ldots, m - 1$ and using the fact that the Euclidean volume of the standard $m - 2$-dimensional simplex is given by $1/(m - 2)!$, we deduce
		\begin{multline}\label{eq_vol_fin_2}
			\int_{t_1 = 0}^{+ \infty} \cdots \int_{t_{m - 1} = 0}^{+ \infty} 
			\int_{x \in \real} \big( x - \frac{1}{k} \sum_{j = 1}^{m - 1} t_j  \overline{\deg}(L_j) \big)_+ d \eta_{k, l}^{HN|C}(x) dt_1 \cdots dt_{m - 1}
			\\
			=
			\frac{k^{m - 1}}{(m - 2)! \cdot \overline{\deg}(L_1) \cdots \overline{\deg}(L_{m - 1})}
			\int_{t = 0}^{+ \infty}
			\int_{x \in \real} ( x - t)_+ t^{m - 2} d \eta_{k, l}^{HN|C}(x) dt.
		\end{multline}
		Tonelli's theorem on another hand yields
		\begin{equation}\label{eq_vol_fin_3}
			\int_{t = 0}^{+ \infty}
			\int_{x \in \real} ( x - t)_+ t^{m - 2} d \eta_{k, l}^{HN|C}(x) dt
			=
			\frac{1}{(m - 1) m}
			\int_{x = 0}^{+ \infty} x^m d \eta_{k, l}^{HN|C}(x).
		\end{equation}
		Remark that we have $\overline{\deg}(L_j) = a_j \cdot l_j$, and so a combination of (\ref{eq_vol_fin_1}), (\ref{eq_vol_fin_2}) and (\ref{eq_vol_fin_3}) yields (\ref{eq_vol_bnd_new}), finishing the proof of (\ref{eq_vol_bnd}).
		\par 
		Assume now that we have $\esssup \eta^{HN} \leq 0$. 
		Then if we have an equality in (\ref{eq_kappa_ineq}), the formula (\ref{eq_vol_bnd}) applies.
		But it implies that ${\rm{vol}}_{\kappa}(L) = 0$, contradicting the definition of the multiplicity. 
		Hence, we have a strict inequality in (\ref{eq_kappa_ineq}).
		\par 
		Finally, let us assume that $\esssup \eta^{HN} < 0$.
		The by the superadditivity from Proposition \ref{prop_superadd}, we deduce that the maximal slopes of $\mathscr{E}_k$, $k \in \nat^*$ are negative.
		From this and Lemma \ref{lem_slope_zero_morph}, for any $k \in \nat^*$, we have $H^0(B, \mathscr{E}_k) = \{0\}$.
		By (\ref{eq_projjj}), we deduce that $L$ is not effective. 
	\end{proof}
	We finish this section by an example, showing that when $\dim B \geq 2$, the equality might not hold in Theorem \ref{thm_vol_bnd}, in contrast with the case $\dim B = 1$.
	\par 
	\textbf{Example}. Consider an elliptic curve $C$ endowed with a non-torsion line bundle $L_0$ of degree $0$.
	Then we clearly have $H^0(C, L_0^{k}) = \{ 0 \}$ for any $k \in \nat^*$.
	Let us now consider the manifold $B := C \times \mathbb{P}^1$, a vector bundle $F := \mathscr{O}_{\mathbb{P}^1}(1) \otimes \comp^2$ over $\mathbb{P}^1$, which we view as a vector bundle over $B$ by the pull-back.
	Let $X := \mathbb{P}(F)$ and define $L$ on $X$ by the tensor product of the relative hyperplane bundle $\mathscr{O}_{\mathbb{P}(F)}(1)$ on $\mathbb{P}(F)$ and the pull-back of $L_0$.
	Then it is immediate that $L$ is relatively ample with respect to the natural submersion $\pi : X \to B$.
	Immediate calculation shows that for any $k \in \nat$, the direct images $\mathscr{E}_k := R^0 \pi_* L^k$ are isomorphic to the tensor product of the pull-backs of $L_0^k$ and ${\rm{Sym}}^k F$.
	In particular, by the Künneth formula, we have $H^0(B, \mathscr{E}_k) = \{ 0\}$, which implies by the projection formula that $H^0(X, L^{k})  = \{ 0\}$, showing in particular that ${\rm{vol}}(L) = 0$.
	\par 
	However, with respect to any multipolarization on $B$, the vector bundles $\mathscr{E}_k$ are semistable (being a product of a line bundle and a trivial vector bundle), and their slopes grow linearly in $k$ with a strictly positive exponent.
	We hence see that the left-hand side of (\ref{eq_vol_bnd}) vanishes, but the right-hand side is positive, and so the equality does not hold in Theorem \ref{thm_vol_bnd}.
	
	\section{Numerical expressions for the asymptotic slopes}\label{sect_num_char}
	The main goal of this section is to give an algebraic interpretation of $\eta_{\min}^{HN}$ and $\eta_{\max}^{HN}$.
	Throughout the whole section, we assume that $L$ is relatively ample.
	\par 
	We say that $\alpha \in H^{1, 1}(X)$ is \textit{$[\omega_B]$-generally fibered nef (resp.  psef) with respect to $\pi$} if there is $l_0 \in \nat^*$, such that for any regular curve $C = C(l) \subset B$, $l = (l_1, \ldots, l_{m - 1})$, $l_i \geq l_0$, defined as before Theorem \ref{thm_mehta_ramanathan}, the restriction of $\alpha$ to $\pi^{-1}(C)$ is nef (resp.  psef), i.e. for any irreducible curve $C'$ in $\pi^{-1}(C)$, the pairing $\int_{C'} \alpha$ is non negative (resp.  the restriction of $\alpha$ to $\pi^{-1}(C)$ contains a positive current).
	For families given by projectivization of vector bundles, an equivalent definition was given by Miyaoka \cite{MiyGenNef}, see also Peternell \cite{PetGenNef}.
	In analogy with the definition of the nef (resp.  psef) threshold, for a class $\alpha \in H^{1, 1}(X)$, we then define the \textit{$[\omega_B]$-general nef (resp.  psef) threshold with respect to $\pi$} as the supremum over all $T \in \real$, such that $\alpha - T (\pi^* [\omega_{B, 0}] / \int_B [\omega_{B, 0}] \cdots [\omega_{B, m - 1}] )$ is $[\omega_B]$-generally fibered nef (resp.  psef) with respect to $\pi$ for some Kähler class $[\omega_{B, 0}] \in H^{1, 1}(B)$ (it is easy to see that the definition is independent on the choice of $[\omega_{B, 0}]$, and for relatively ample $\alpha$, this supremum is finite).	
	\par 
	We can now state the main result of this section.
	\begin{thm}\label{cor_treshhold}
		Assume $L$ is relatively ample and the multipolarization $[\omega_B]$ is integral.
		Then the quantity $\eta_{\min}^{HN}$ (resp. $\eta_{\max}^{HN}$) equals to $[\omega_B]$-general nef (resp. psef) threshold of $c_1(L)$ with respect to $\pi$.
	\end{thm}
	\par 
	In order to prove Theorem \ref{cor_treshhold}, let us first recall the results of Xu-Zhuang from \cite{XuZhuPosCM}.
	We conserve the notations from the introduction, and assume for the moment that $\dim B = 1$.
	For a class $\alpha \in H^{1, 1}(X)$, we then define the \textit{nef (resp.  psef) threshold of $\alpha$ with respect to $\pi$} as the supremum over all $T \in \real$, such that $\alpha - T (\pi^* [\omega_{B, 0}] / \int_B [\omega_{B, 0}] )$ is nef (resp.  psef) for some Kähler class $[\omega_{B, 0}]$ on $B$.	
	\begin{prop}\label{prop_xu_zh}
		If $\dim B = 1$, then $\eta_{\min}^{HN}$ (resp. $\eta_{\max}^{HN}$) equals to nef (resp. psef) threshold of $c_1(L)$ with respect to $\pi$.
		Moreover, we have $\eta_{\min}^{HN} = \essinf \eta^{HN}$.
	\end{prop}
	\begin{proof}
		In \cite[Proposition 2.28]{XuZhuPosCM}, authors proved that $\esssup \eta^{HN}$ equals to psef threshold of $c_1(L)$ with respect to $\pi$.
		From this and Theorem \ref{thm_filt}, we deduce the statement of Proposition \ref{prop_xu_zh} for $\eta_{\max}^{HN}$.
		In the proof of \cite[Lemma 2.26]{XuZhuPosCM}, authors established that $\eta_{\min}^{HN}$ is not smaller than nef threshold of $c_1(L)$ with respect to $\pi$.
		This finishes the proof of Proposition \ref{prop_xu_zh}, as by Theorem \ref{thm_filt}, we have $\eta_{\min}^{HN} \leq \essinf \eta^{HN}$, and by \cite[Proposition 2.28]{XuZhuPosCM}, $\essinf \eta^{HN}$ equals to nef threshold of $c_1(L)$ with respect to $\pi$.
	\end{proof}
	\par
	Let us now deduce the following numerical formulas
	\begin{equation}\label{eq_num_form}
		\eta_{\min}^{HN}
		=
		\inf_{C \in \mathcal{C}_{irr}} \frac{\int_C c_1(L)}{{\deg(\pi|_C)}}
		,
		\qquad
		\eta_{\max}^{HN}
		=
		\inf_{C \in \mathcal{C}_{mov}} \frac{\int_C c_1(L)}{{\deg(\pi|_C)}},
	\end{equation}
	where $\mathcal{C}_{irr}$ is the set of irreducible curves $C \subset X$, such that $\pi|_C$ is surjective, $\mathcal{C}_{mov} \subset \mathcal{C}_{irr}$ is the subset of movable curves, i.e. curves which can be put in a family, which cover $X$, and $\deg(\pi|_C)$ is the topological degree of $\pi|_C$.
	 By the definition of nefness, we see that the nef threshold of $c_1(L)$ with respect to $\pi$ equals
	\begin{equation}
		\inf_{C \in \mathcal{C}_{irr}} \frac{\int_C c_1(L) }{\int_C \pi|_C^* [\omega_{B, 0}] / \int_B [\omega_{B, 0}]},
	\end{equation}
	where $\mathcal{C}_{irr}$ was defined before (\ref{eq_num_form}).
	This gives us the first formula from (\ref{eq_num_form}) by the definition of the topological degree.
	The proof of the second formula is identical, except that it is necessary to use the characterization of the dual to the pseudoeffective cone in terms of movable curves due to Boucksom-Demailly-Păun-Peternell \cite[Theorem 2.2]{BDPP}.
	\begin{sloppypar}
	\begin{proof}[Proof of Theorem \ref{cor_treshhold}]
		Directly from Proposition \ref{prop_xu_zh}, we have 
		 \begin{equation}\label{eq_min_hn00}
		 	\eta_{\min, l}^{HN|C}
		 	=
		 	\sup
		 	\bigg\{
		 	T \in \real : \text{ the class } 
		 	c_1(L)|_{\pi^{-1}(C)} - T \cdot \frac{ \pi^* [\omega_{B, 0}]|_{\pi^{-1}(C)} }{ \int_C [\omega_{B, 0}]} \text{ is nef}
		 	\bigg\},
		 \end{equation}
		 where the curves $C = C(l)$ are as in Theorem \ref{thm_mehta_ramanathan}.
		 Similar formulas hold for the maximal asymptotic slope.
		 We deduce Theorem \ref{cor_treshhold} by (\ref{eq_min_hn00}), the identity $\int_C [\omega_{B, 0}] = l_1 \cdots l_{m - 1} \cdot \int_B [\omega_{B, 0}] \cdots [\omega_{B, m - 1}]$ and Theorem \ref{thm_mehta_ramanathan}.
	\end{proof}
	\end{sloppypar}
	Let us now give another interpretation of Theorem \ref{cor_treshhold}.
	We say that a relatively ample $\mathbb{Q}$-line bundle $L$ on $X$ is \textit{stably $([\omega_{B, 1}], \ldots, [\omega_{B, m - 1}])$-generally fibered nef with respect to $\pi$} if for some (or any) ample line bundle $L_0$ on $X$, for any $\epsilon > 0$, $\epsilon \in \mathbb{Q}$, the $\mathbb{Q}$-line bundle $L \otimes L_0^{\epsilon}$ is $([\omega_{B, 1}], \ldots, [\omega_{B, m - 1}])$-generally fibered nef with respect to $\pi$.
	\begin{prop}\label{prop_interpr_st_gen_pos}
		A relatively ample line bundle $L$ is stably $[\omega_B]$-generally fibered nef with respect to $\pi$ if and only if $\eta_{\min}^{HN} \geq 0$.
	\end{prop}
	\begin{proof}
		We assume first that $\eta_{\min}^{HN} \geq 0$.
		By Theorem \ref{cor_treshhold}, for any $\epsilon > 0$, $\epsilon \in \mathbb{Q}$, the $\mathbb{Q}$-line bundle $L \otimes \pi^* L_1^{\epsilon}$ is $[\omega_B]$-generally fibered nef with respect to $\pi$, where $L_1$ is an ample line bundle on $B$.
		But then for any ample line bundle $L_0$ on $X$, any $\epsilon > 0$, $\epsilon \in \mathbb{Q}$, the $\mathbb{Q}$-line bundle $L \otimes L_0^{\epsilon}$ is $[\omega_B]$-generally fibered nef with respect to $\pi$, which means that $L$ is stably $[\omega_B]$-generally fibered nef with respect to $\pi$.
		\par 
		Inversely, assume that $L$ is stably $[\omega_B]$-generally fibered nef with respect to $\pi$.
		Let $N \in \nat^*$ be such that the line bundle $L_0 := L \otimes \pi^* L_1^N$ is ample (such $N$ exists since $L$ is relatively ample).
		Then for any $\epsilon > 0$, $\epsilon \in \mathbb{Q}$, the $\mathbb{Q}$-line bundle $L \otimes L_0^{\epsilon}$ is $[\omega_B]$-generally fibered nef with respect to $\pi$.
		However, since $L \otimes L_0^{\epsilon} = L^{1 + \epsilon} \otimes \pi^* L_1^{N \epsilon}$, we see that the $[\omega_B]$-general fibered nef threshold of $c_1(L)$ is non-negative, which by Theorem \ref{cor_treshhold} means that $\eta_{\min}^{HN} \geq 0$, and this finishes the proof.
	\end{proof}
	\begin{sloppypar}
	\begin{rem}
		From the proof, we see that in the definition of stably $[\omega_B]$-generally fibered nefness with respect to $\pi$, instead of ample line bundle $L_0$ over $X$, we can consider the pull-back of an ample line bundle over the base.
	\end{rem}
	\end{sloppypar}

\section{Possible directions for further research}\label{sect_open}
	In this section we will discuss some open problems concerning the relation between the current work and the valuations theorems \cite{ReesValI}, \cite{ReesValII} and \cite{BouckJohn21}.
	We will assume now that $L$ is relatively ample.
	We conserve the notations from Introduction.
	We pick $b \in B$ as described before (\ref{eq_ind_measu_defn_HN}), and denote by $\mathcal{F}_b^{HN}$ the filtration on $R(X_b, L_b)$ induced by the Harder-Narasimhan filtration on $\mathscr{E}_k$ as in Section \ref{sect_sm_bnd}.
	By Proposition \ref{prop_sm}, the resulting filtration is submultiplicative. 
	\par 
	We denote by $X_b^{{\textrm{an}}}$ the Berkovich analytification of $X_b$, which is the set of all semivaluations on $X_b$, i.e. $\mathbb{R}$-valued valuations on the function field of all subvarieties of $X_b$, trivial on $\comp$. 
	Following \cite{BouckJohn21}, we define a function $FS^{HN}_{k, b}: X_b^{{\textrm{an}}} \to \mathbb{R}$, as $v \mapsto \max_{s \in H^0(X_b, L_b^{k})} (- v(s) + w_{\mathcal{F}_{k, b}^{HN}}(s))$, where $w_{\mathcal{F}_{k, b}^{HN}}$ is the weight of the filtration.
	By the submultiplicativity of $\mathcal{F}_b^{HN}$, the sequence of functions $FS^{HN}_{k, b}$ is superadditive, cf. \cite{BouckJohn21}, and so, as $k \to \infty$, the sequence $\frac{1}{k} FS^{HN}_{k, b}$ converges pointwise to a function, which we denote by $FS^{HN}_{b}: X_b^{{\textrm{an}}} \to \mathbb{R} \cup \{+ \infty \}$, and call the Fubini-Study potential of the Harder-Narasimhan filtration at $b \in B$.
	By the boundness of $\mathcal{F}_b^{HN}$, cf. Proposition \ref{thm_bound}, the function $FS^{HN}_{b}$ is bounded (and so the value $+ \infty$ is not attained). 
	\par 
	\textbf{Question 1}: Is the function $FS^{HN}_{b}$ continuous (when the space $X_b^{{\textrm{an}}}$ is endowed with the topology of pointwise convergence)?
	\par
	It seems that addressing the above question is important for understanding finer the asymptotic properties of the Harder-Narasimhan filtration.
	To explain this, let us recall some general theory of submultiplicative filtrations.
	We fix a complex projective manifold $Y$ with an ample line bundle $L$.
	Recall, cf. \cite[(2.11)]{BouckJohn21}, that an arbitrary bounded function $\phi : Y^{{\textrm{an}}} \to \mathbb{R}$ induces the filtration ${\rm{IN}}(\phi) = \oplus_{k = 0}^{+\infty} {\rm{IN}}(\phi)_k$ on $R(Y, L)$, defined so that the weight function $w_{{\rm{IN}}(\phi)_k} : H^0(Y, L^{k}) \to \real$ is given by
	\begin{equation}
		w_{{\rm{IN}}(\phi)_k}(s) := \inf_{v \in Y^{{\textrm{an}}}} \{ v(s) + k \phi(v) \}. 
	\end{equation}
	An easy verification shows that the function $w_{{\rm{IN}}(\phi)_k}$ corresponds to a decreasing filtration on $H^0(Y, L^{k})$; this simply means that $w_{{\rm{IN}}(\phi)_k}(\lambda s) = w_{{\rm{IN}}(\phi)_k}(s)$ for $\lambda \neq 0$, and $w_{{\rm{IN}}(\phi)_k}(s_1 + s_2) \geq \min \{ w_{{\rm{IN}}(\phi)_k}(s_1), w_{{\rm{IN}}(\phi)_k}(s_2) \}$.
	These filtrations altogether form a filtration on $R(Y, L)$, which we denote by ${\rm{IN}}(\phi)$, following the notations from \cite[\S 2.4]{BouckJohn21}.
	An easy verification shows that ${\rm{IN}}(\phi)$ is a submultiplicative filtration.
	The asymptotic properties of submultiplicative filtrations of the form ${\rm{IN}}(\phi)$ has been studied extensively in the past, cf. \cite[Theorem 6.4.6]{ChenMoriwAdelic}, \cite{BouckJohn21}.
	\par 
	It is then natural to compare, for a given submultiplicative filtration $\mathcal{F}$ on $R(Y, L)$, the filtrations $\mathcal{F}$ and ${\rm{IN}}(FS_{\mathcal{F}})$. 
	The valuation theorems of Rees \cite{ReesValI}, \cite{ReesValII} tell that if the filtration $\mathcal{F}$ is of finite type, then there is $C > 0$ such that for any $k \in \nat$, the difference between the weight functions on the two filtrations, restricted to $H^0(Y, L^{k})$, is bounded by $C$.
	Now, for filtrations of finite type, $FS_{\mathcal{F}}$ is continuous, cf. \cite[Theorem 2.16]{BouckJohn21}. Boucksom-Jonsson in \cite[Theorems 2.16, 2.19 and Proposition 2.24]{BouckJohn21} more generally established that if $FS_{\mathcal{F}}$ is continuous, then for any $\epsilon > 0$, there is $k_0 \in \nat$, so that the difference between the weight functions on the two filtrations, restricted to $H^0(Y, L^{k})$, is bounded by $\epsilon k$ for any $k \geq k_0$.
	If $\mathcal{F}$ is only bounded, the last comparison fails, cf. \cite[Example 2.25]{BouckJohn21}.
	However, in \cite[Theorem 3.18]{BouckJohn21}, authors established that boundness suffices to establish that the two filtrations $\mathcal{F}$ and ${\rm{IN}}(FS_{\mathcal{F}})$ have the same asymptotic jumping measure.
	In this vein, the following question might be seen as a far-reaching refinement of Question 1.
	\par 
	\textbf{Question 2}: Is the filtration $\mathcal{F}_b^{HN}$ finitely generated?
	\par 
	It seems that due to irrationality issues it is obligatory to assume in the above question that $[\omega_B] \in H^2(B, \mathbb{Z})$.
	Continuing in this direction, if the answer to Question 2 is positive, it would be very interesting to give an analytic description of the test configuration associated with $\mathcal{F}_b^{HN}$ (recall that finitely generated filtration with integer weights are more or less in one-to-one correspondence with test configurations, cf. \cite{NystOkounTest}, \cite{SzekeTestConf}, \cite{BouckHisJonsFourier}). 
	Probably, it will be related with the Wess-Zumino-Witten equation, as the variation of the weights of the filtration $\mathcal{F}_b^{HN}$ provide obstructions to the solutions of it, see \cite{FinHNII}, \cite{FinHYM}.
	\par 
	Now, it is very tempting to believe that one can approach the questions of the current article using the above non-Archimedean framework.
	To explain this more precisely, let us introduce some further notations.
	By our choice of $b \in B$, prompted by Proposition \ref{prop_exist_b}, for any $l \in \nat^{*(m - 1)}$, we can pick a curve $\iota_l: C_l \hookrightarrow B$, passing through $b \in B$, such that for any $k \in \nat$, the Harder-Narasimhan slopes of the pull-back of $\mathscr{E}_k$ to $C_l$ coincide with the Harder-Narasimhan slopes of the pull-back of $\mathscr{E}_k$ to a general curve given by a complete intersection in $l [\omega_B]$.
	We denote by $\mathcal{F}_{k, b, l}^{HN}$ the filtration induced by the Harder-Narasimhan filtration on $\iota_l^* \mathscr{E}_k$ on the fiber $H^0(X_b, L_b^{k})$ of $\mathscr{E}_k$ at $b \in C_l$.
	For a fixed $l$, these filtrations form a submultiplicative filtration on $R(X_b, L_b)$, which we denote by $\mathcal{F}_{b, l}^{HN} = \oplus_{k = 0}^{\infty} \mathcal{F}_{k, b, l}^{HN}$.
	We denote by $FS^{HN}_{b, l}: X_b^{{\textrm{an}}} \to \mathbb{R} \cup \{+ \infty \}$, the Fubini-Study potential associated with the filtration given by the division of the slopes of $\mathcal{F}_{b, l}^{HN}$ by $l_1 \cdots l_{m - 1}$.
	\par 
	\textbf{Question 3}: Does the sequence of potentials $FS^{HN}_{b, l}$ converge uniformly to $FS^{HN}_{b}$, as $l \to \infty$?
	If not uniformly, does the convergence hold in $d_p$-topology, $p \in [1, +\infty[$, introduced in \cite[\S 5]{BouckJohn21}?
	\par 
	If the answer to Questions 1 and 3 are positive, then by relying on the techniques of \cite{BouckJohn21}, one can give an alternative proof of Theorem \ref{thm_mehta_ramanathan}.
	Remark also that the maximal asymptotic slope of filtrations of the form ${\rm{IN}}(\phi)$, for $\phi$ continuous, admit a description in terms of successive minima, see \cite[Proposition 8.10.1]{ChenMoriwHilSam}, which are defined in terms of the quotient filtrations associated with subvarieties.
	If the answer to Question 1 is positive, the study of successive minima would apply to Harder-Narasimhan filtration by \cite{BouckJohn21}. 
	It might be interesting to explore the connection between this description of the maximal asymptotic slope and the one given by Theorem \ref{cor_treshhold}.
	\par 
	Let us, however, point out a possible obstacle in addressing Question 3.
	While Corollary \ref{cor_welldef_slo} tells us that the jumping numbers of  slopes of $\mathcal{F}_{b, l}^{HN}$ do not depend on the choice of the curve $C_l$ (which makes the measure $\eta_l^{HN|C}$ from Theorem \ref{thm_mehta_ramanathan} well-defined), it says nothing on the dependence of the filtration itself.
	It is, hence, a priori not clear if the potentials $FS^{HN}_{b, l}$ depend on the choice of $C_l$, which raises doubts concerning the possible strengthening on the level of Fubini-Study potentials of the subadditivity property established in Theorem \ref{thm_coord_subm}. 
	
\bibliography{bibliography}

		\bibliographystyle{abbrv}

\Addresses

\end{document}